\documentclass[a4paper]{amsart}
\usepackage[utf8]{inputenc}
\usepackage{tikz}
\usepackage{url}
\usetikzlibrary{decorations.markings}
\usetikzlibrary{positioning}

\title{On a vertex-minimal triangulation of $\mathbb{R}P^4$}
\date{}
\author{Sonia Balagopalan}
\address{Department of Mathematics and Statistics, National University of Ireland Maynooth, Maynooth, Co. Kildare, Ireland}
\email{s.balagopalan@gmail.com}
\curraddr{Institute of Mathematics, the Hebrew University of Jerusalem, Givat Ram, Jerusalem, 91904, Israel.}
\newcommand{\RR}{\mathbb{R}}
\newcommand{\ee}{\mathbf{e}}

\begin{document}

\begin{abstract}
 We give three constructions of a vertex-minimal triangulation of $4$-dimensional real projective space 
 $\mathbb{R}P^4$. The first construction describes a $4$-dimensional sphere on $32$ vertices, which is 
 a double cover of a triangulated $\mathbb{R}P^4$ and has a large amount of symmetry. The second and third 
 constructions illustrate  approaches to improving the known number of vertices needed to triangulate 
 $n$-dimensional real projective space. All three constructions deliver the same combinatorial manifold, 
 which is also the same as the only known $16$-vertex triangulation of $\mathbb{R}P^4$. We also give a short, 
 simple construction of the $22$-point Witt design, which is closely related to the complex we construct.
\end{abstract}


\maketitle

\newtheorem{thm}{Theorem}
\newtheorem*{conj}{Conjecture}
\theoremstyle{remark}	  \newtheorem{ex}{Example}
\theoremstyle{remark}	  \newtheorem*{rk}{Remark}
\theoremstyle{definition} \newtheorem{cons}{Construction}

\section{Introduction}

How many vertices does it take to (simplicially) triangulate real projective $n$-space, $\RR P^n$? 
It is well known that the answer is $6$ when $n=2$, with the triangulation realized as the antipodal 
quotient of the icosahedron. In 1969, D.W. Walkup proved that a vertex-minimal triangulation of $\RR P^3$ 
requires $11$ vertices, and described all such triangulations~\cite{W}. P. Arnoux and A. Marin, in 1991, 
proved that the minimum number of vertices needed to triangulate $\mathbb{R}P^n,\ n\geq 3$, is 
at least $\binom{n+2}{2}+1$\cite{AM}.

In 1987, W. K\"uhnel gave a triangulation of $\RR P^n$ using $2^{n+1}-1$ vertices, which takes the 
barycentric subdivision of the boundary of the $n+1$-simplex and quotients it by an antipodal map~\cite{K}. 
The K\"uhnel construction gives the smallest known explicit triangulations of $\RR P^n$ for $n>5$. For a 
survey of these and other results on minimal triangulations, and also all relevant definitions, see~\cite{D}.

The {\tt{BISTELLAR}} program of F.H. Lutz uses a heuristic search algorithm to reduce the $f$-vector 
of a given complex using bistellar flips~\cite{BL,L}. Among the several combinatorial manifolds found by 
this program was a $16$-vertex triangulation of $\mathbb{R}P^4$, called $\RR P^4_{16}$~\cite[p.77]{L}. This
complex was obtained by applying the {\tt{BISTELLAR}} program to the $31$-vertex $\RR P^4$ due to K\"uhnel. 
The automorphism group of this particular triangulation was also computed in~\cite[p.77]{L}, and was found 
to be $S_6$, which acts on the $16$-element vertex set by splitting it into orbits of size $6$ and $10$, and 
on the set of $150$ facets by splitting it into orbits of size $30$ and $120$. No other triangulation of 
$\RR P^4$ on $16$ vertices is known. Apart from the information described above, there does not seem 
to be anything else known about a $16$-vertex $\RR P^4$ in the literature. We believe that this extremal object 
exhibiting such a high degree of symmetry deserves to be better understood, both for its own sake and for furthering
our understanding of the general problem of triangulating $\RR P^n$.

We present three constructions of a triangulated $\RR P^4$ on ${16}$ vertices, all of which turn out to be 
isomorphic to $\RR P^4_{16}$, and make some observations about the remarkable combinatorial structure of this complex.

Construction~\ref{c1} describes $\RR P^4_{16}$ as the quotient of a triangulated 
$4$-sphere on $32$ vertices. In order for such a construction to be possible, the $S^4$ we 
construct needs to be antipodal. We say that a simplicial complex $K$ is \emph{antipodal} if it
is invariant under an involution $\sigma$ such that the (graph) distance between vertices $v$ and $\sigma(v)$ 
in the $1$-skeleton of $K$ is at least $3$. In particular, the links of $v$ and $\sigma(v)$ are disjoint,
and isomorphic to the link of $\overline{v}$ in the quotient complex $K/\langle\sigma\rangle$, and we say that $\sigma$ is 
\emph{link-separating} on $k$. If $K$ is a combinatorial manifold, then it follows that $K/\langle\sigma\rangle$ is also a
combinatorial manifold. We call $\sigma$ the \emph{antipodal map}. To be precise, when we refer to an antipodal complex, we 
are implicitly refering to a pair $(K,\sigma)$.

We also point out a connection between the triangulated $\RR P^4$ of Construction~\ref{c1}, and the 
Witt design on $22$ points, $W_{22}$. We note that the smaller orbit of the facets of $\RR P^4_{16}$ under its 
automorphism group is closely related to a well known symmetric $2-$design or \emph{biplane} on $16$ points. Also, 
the partition of the vertex set into orbits suggests the construction of a ``dual'' biplane by the introduction of six new 
points. The resulting configuration can be extended to a $3-$design on $22$ points, with $77$ blocks of size $6$, of which 
$W_{22}$ is unique up to isomorphism. 

Our construction of $W_{22}$ is short and elementary, and does not seem to appear as such in the literature. Nevertheless, 
no construction of an object so well understood can be said to be entirely new, and ours has many factors in common to two 
previously known contstructions, which we briefly note. We justify our choice to present our construction in full detail, 
more for what it illuminates about $\RR P^4_{16}$ and its automorphism group, than for what it says about $W_{22}$.

Construction~\ref{c1}, though it exposits the remarkable symmetries of  $\RR P^4_{16}$, seems to rely on exceptional properties
of the number $6$, and is not a very encouraging as a model for analogous constructions in higher dimensions. We give two more 
constructions which are more promising in this direction. These constructions view $\RR P^4$ as a $4$-dimensional ball with 
antipodal simplices on its boundary identified. We start with a suitable convex $4$-polytope, and place it inside its dual, 
and triangulate the regions thus formed, till we get a $3$-sphere on the boundary, which we can glue to itself to give an 
$\RR P^4$. Our second and third constructions follow this strategy, starting from a $16$-cell and a suspended cube respectively. 
We also describe a way of looking at Walkup's $\RR P^3_{11}$ and even $\RR P^2_{6}$ as $3$ and $2$-dimensional analogues of our 
constructions. Both these constructions give the same complex as the first. A simple observation about the automorphism groups 
of these complexes allows us to identify all three constructions with each other. This supports the conjectured uniqueness of 
$\RR P^4_{16}$ as the vertex-minimal triangulation of $\RR P^4$. 

\section{First construction and combinatorial properties}

We construct $\RR P^4_{16}$ by starting with the standard $4$-sphere and constructing an antipodal 
$S^4$ on $32$ vertices, by successive transformations.

\begin{cons} \label{c1} 
Let $\Delta^5$ denote the standard $5$-simplex in $\RR^6$. Let $\ee_i, 1\leq i \leq 6$ denote 
the $i^{\textrm{th}}$ elementary vector. The set $V_1=\{\ee_i|1\leq i\leq 6\}$ is the vertex set 
of $\Delta^5$.  The boundary of $\Delta^5$ is a triangulated $4$-sphere on these six vertices, and 
each of its facets is a $4$-simplex containing five elements of $V_1$. Call this complex $X_{6}$.

Let $\mathbf{1}=\sum_{i=1}^6\ee_i$. Then the barycenter of the facet $\Delta_i$ with vertex set 
$V_1\setminus \ee_i$ of $X_6$ is the point $\frac{1}{5}(\mathbf{1}-\ee_i)\in \RR^6$, where 
$1\leq i\leq 6$. Call the set of these points $V_5$. Introducing these points allows us to subdivide 
each facet of $X_6$ as the union of five $4$-simplices, by replacing $\Delta_i$ with the cone over 
each of its tetrahedra at the point $\frac{1}{5}(\mathbf{1}-\ee_i)$. This gives a 12-vertex triangulated 
$S^4$ with vertex set $V_1\cup V_5$. The facets of this complex are all the $4$-simplices of the form 
$[\ee_i,\ee_j,\ee_k,\ee_l,\frac{1}{5}(\mathbf 1-\ee_m)]$, where $i,j,k,l,m$ are distinct. Call this complex 
$X_{12}$.

Let $V_3=\{\frac{1}{3}(\ee_i+\ee_j+\ee_k)|1\leq i<j<k\leq6\}$ denote the set of barycenters of the triangles 
of $\Delta^5$. We use these points to further subdivide each facet of $X^4_{12}$ in the following way. The 
tetrahedron $[\ee_i,\ee_j,\ee_k,\ee_l]$ of the facet $[\ee_i,\ee_j,\ee_k,\ee_l,\frac{1}{5}(\mathbf{1}-\ee_m)]$ 
can be decomposed into eleven tetrahedra. These are, from the outside in, six tetrahedra of the form 
$[\ee_i,\ee_j,\frac{1}{3}(\ee_i+\ee_j+\ee_k),\frac{1}{3}(\ee_i+\ee_j+\ee_l)]$
corresponding to every pair of elements of $\{i,j,k,l\}$, four tetrahedra of the form 
$[\ee_i,\frac{1}{3}(\ee_i+\ee_j+\ee_k),\frac{1}{3}(\ee_i+\ee_j+\ee_l),\frac{1}{3}(\ee_i+\ee_k+\ee_l)]$ 
corresponding to every element of $\{i,j,k,l\}$, and the tetrahedron
$[\frac{1}{3}(\ee_i+\ee_j+\ee_k),\frac{1}{3}(\ee_i+\ee_j+\ee_l),\frac{1}{3}
(\ee_i+\ee_k+\ee_l),\frac{1}{3}(\ee_j+\ee_k+\ee_l)]$. See Figure~\ref{subtet} for an illustration of this subdivision. 
Here, $\ee_{ijk}$ denotes $\frac{1}{3}(\ee_i+\ee_j+\ee_k)\in V_3$.

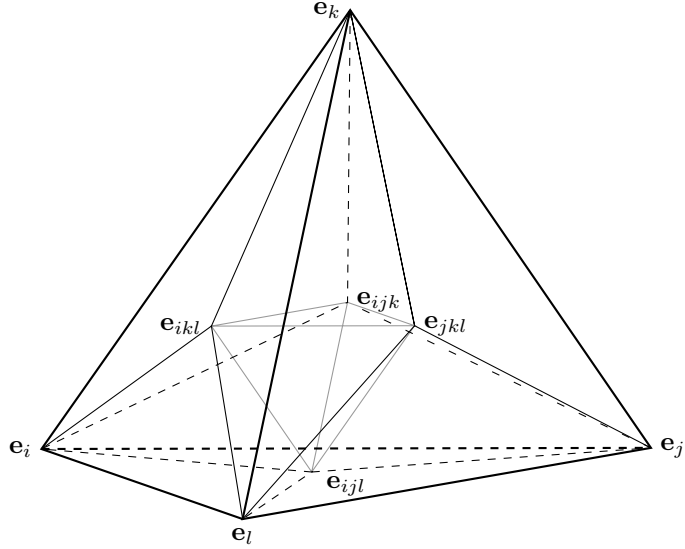
\begin{figure}
\begin{center}
\begin{tikzpicture}[z=.3,y=24]
\coordinate [label=left:$\ee_i$] (i) at (-4,-2.828,-2.828)  ;  
\coordinate [label=right:$\ee_j$] (j) at (4,-2.828,-1.414) ; 
\coordinate [label=left:$\ee_k$] (k) at (0,4,2.828) ; 
\coordinate [label=below:$\ee_l$] (l) at (-1.414,-4,2.828) ;
\coordinate [label={[xshift=2pt,yshift=2pt]below right:$\ee_{ijl}$}] (ijl) at (-0.471, -3.219, -0.471) ;  
\coordinate [label=right:$\ee_{ijk}$] (ijk) at (0, -0.552, -0.471); 
\coordinate [label=left:$\ee_{ikl}$] (ikl) at (-1.805, -0.943, 0.943);
\coordinate [label=right:$\ee_{jkl}$] (jkl) at (0.862, -0.943, 1.414);

\draw [black!40](ijk)-- (jkl)--(ikl)--(ijl)--(jkl) (ijl)--(ijk)--(ikl);
\draw (j)--(jkl)--(k) (jkl)--(l)--(ikl)--(k)--(jkl) (i)--(ikl);
\draw [dashed](k)--(ijk)--(i)--(ijl)--(l) (ijl)--(j);
\draw [loosely dashed](j)--(ijk);
\draw [thick,dashed](i)--(j);
\draw [thick] (k)--(l)--(j)--(k)--(i)--(l) ;

\end{tikzpicture}
\caption{Tetrahedron subdivided by barycenters of its $2$-faces}
\label{subtet}
\end{center}
\end{figure}

We take the join of $\frac{1}{5}(\mathbf{1}-\ee_m)$ with each of these tetrahedra to obtain a decomposition of the facet.
This gives us a triangulated $S^4$ on $32$ vertices, $X_{32}$, with three kinds of facets, containing two, three, and 
four vertices of $V_3$ respectively. The vertex set $V_1\cup V_3\cup V_5$ of $X_{32}$ suggests a natural choice 
of antipodal map, the one that swaps $\ee_i$ with $\frac{1}{5}(\mathbf{1} -\ee_i)$ and $\frac{1}{3}(\ee_i+\ee_j+\ee_k)$ 
with $\frac{1}{3}(\mathbf{1}-\ee_i-\ee_j-\ee_k)$.

In order to obtain an antipodal complex from  $X_{32}$,  we need to transform the complex to one that is invariant 
under the above map, and also separate antipodal vertices, till they are far enough apart.  We do this using bistellar 
flips. 

In $X_{32}$, any two vertices of $V_1$ form an edge, but no two vertices of $V_5$ do. Also, any vertex in $V_1$ is 
adjacent to any vertex in $V_5$ other than its antipode. By separating each of the edges within $V_1$, we can reduce 
the asymmetry of the complex, and also increase the distance between would-be antipodal pairs in $V_1\cup V_5$. To achieve 
this, first note that each edge $[\ee_i,\ee_j]$ is contained in four triangles of the form $[\ee_i,\ee_j,
\frac{1}{5}(\mathbf{1}-\ee_k)]$, corresponding to each of its four neighbouring facets $\Delta_k$ in $X_6$. The link 
of $[\ee_i,\ee_j,\frac{1}{5}(\mathbf{1}-\ee_k)]$ in $X_{32}$ is the boundary of the triangle 
$[\frac{1}{3}(\ee_i+\ee_j+\ee_{l}), \frac{1}{3}(\ee_i+\ee_j+\ee_{l'}), \frac{1}{3} (\ee_i+\ee_j+\ee_{l''})]$, where
$\{i,j,k,l,l',l''\}=\{1,\ldots,6\}$, each edge of said triangle corresponding to a tetrahedron containing $[\ee_i,\ee_j]$
joined with $\frac{1}{5}(\mathbf{1}-\ee_k)$ in $X_{12}$. Since any triangle in $X_{32}$ with vertices from $V_3$ is 
obtained by subdividing a tetrahedron, a triangle of the above type is not a face of $X_{32}$.

On the other hand, every triple of vertices of the form 
$\{\frac{1}{3}(\ee_i+\ee_j+\ee_{l}), \frac{1}{3}(\ee_i+\ee_j+\ee_{l'}), \frac{1}{3}(\ee_i+\ee_j+\ee_{l''})\}$ 
is the vertex set of the link of a unique triangle of the form $[\ee_i,\ee_j,\frac{1}{5}(\mathbf{1}-\ee_k)]$,
since $i,j,k,l,l',l''$ are all distinct. So we can apply simultaneous bistellar flips to all triangles of the form
$[\ee_i,\ee_j,\frac{1}{5}(\mathbf{1}-\ee_k)]$, replacing the facets of the form
$$[\ee_i,\ee_j,\frac{1}{5}(\mathbf{1}-\ee_k)]*\partial[\frac{1}{3}
(\ee_i+\ee_j+\ee_{l}), \frac{1}{3}(\ee_i+\ee_j+\ee_{l'}), \frac{1}{3}
(\ee_i+\ee_j+\ee_{l''})]$$ 
with the facets of the form
$$[\frac{1}{3} (\ee_i+\ee_j+\ee_{l}), \frac{1}{3}(\ee_i+\ee_j+\ee_{l'}),
\frac{1}{3} (\ee_i+\ee_j+\ee_{l''})]*\partial[\ee_i,\ee_j,\frac{1}{5}(\mathbf{1}
-\ee_k)].$$
After this round of flips, the link of $[\ee_i,\ee_j]$ is the boundary of the tetrahedron
$$[\frac{1}{3} (\ee_i+\ee_j+\ee_{l}), \frac{1}{3}(\ee_i+\ee_j+\ee_{l'}),
\frac{1}{3} (\ee_i+\ee_j+\ee_{l''}), \frac{1}{3}(\ee_i+\ee_j+\ee_{l'''})]$$
where $i,j,l,l',l'',l'''$ are distinct, each face of said tetrahedron being one introduced 
in the previous round of flips. As above, we can perform simultaneous bistellar flips to replace 
$$[\ee_i,\ee_j]*\partial[\frac{1}{3}(\ee_i+\ee_j+\ee_{l}),\frac{1}{3} (\ee_i+\ee_j+\ee_{l'}),
\frac{1}{3} (\ee_i+\ee_j+\ee_{l''}), \frac{1}{3}(\ee_i+\ee_j+\ee_{l'''})]$$
with
$$[\frac{1}{3} (\ee_i+\ee_j+\ee_{l}), \frac{1}{3}(\ee_i+\ee_j+\ee_{l'}),
\frac{1}{3} (\ee_i+\ee_j+\ee_{l''}), \frac{1}{3}(\ee_i+\ee_j+\ee_{l'''})]*\partial[\ee_i,\ee_j]$$

The facets of the resulting complex are of four types, (or two types under the action of $C_2\times
S_6$, where $S_6$ permutes the axes of $\RR^6$ and $C_2$ is generated by the antipodal involution): 
$$[\ee_i,\frac{1}{3} (\ee_i+\ee_j+\ee_{l}), \frac{1}{3}(\ee_i+\ee_j+\ee_{l'}),
\frac{1}{3} (\ee_i+\ee_j+\ee_{l''}), \frac{1}{3}(\ee_i+\ee_j+\ee_{l'''})]$$
($6\times 5=30$ in number),
$$[\frac{1}{5}(\mathbf{1}-\ee_m),\frac{1}{3} (\ee_i+\ee_j+\ee_{k}),
\frac{1}{3}(\ee_i+\ee_j+\ee_{l}),
\frac{1}{3} (\ee_i+\ee_k+\ee_{l}), \frac{1}{3}(\ee_j+\ee_k+\ee_{l})]$$
($i,j,k,l,m$ distinct, $30$ in number),
$$[\ee_i,\frac{1}{5}(\mathbf{1}-\ee_k),\frac{1}{3} (\ee_i+\ee_j+\ee_{l}),
\frac{1}{3}(\ee_i+\ee_j+\ee_{l'}),\frac{1}{3}(\ee_i+\ee_j+\ee_{l''})]$$
($i,j,k,l,l',l''$ distinct, $6\times 5\times 4=120$ in number),
$$[\ee_i,\frac{1}{5}(\mathbf{1}-\ee_m),\frac{1}{3} (\ee_i+\ee_j+\ee_{k}),
\frac{1}{3}(\ee_i+\ee_j+\ee_{l}),\frac{1}{3}(\ee_{i}+\ee_k+\ee_{l})]$$
($i,j,k,l,m$ distinct, $6\times 5\times 4=120$ in number). 

This complex, call it $S^4_{32}$, is antipodal under the involution that takes $\ee_i$ to 
$\frac{1}{5}(\mathbf{1}-\ee_i)$ and $\frac{1}{3}(\ee_i+\ee_j+\ee_k)$ to $\frac{1}{3}(\mathbf{1}-\ee_i-\ee_j-\ee_k)$, 
and commutes with the $S_6$ action induced from the permutations of $V_1$. This involution is link-separating on $S^4_{32}$, 
and so quotients it to give a triangulated $\RR P^4$ on $16$ vertices. \qed
\end{cons}

If, in the above construction, we deform the $5$-simplex to choose the elements of $V_3$ and $V_5$ to be of the form 
$e_i+e_j+e_k$ and $\mathbf{1}-e_i$ respectively, the antipodal map is just $\mathbf{x}\mapsto\mathbf{1}-\mathbf{x}$ 
on $V_1\cup V_3\cup V_5$, which also acts on the analogous geometric carrier of the complex above. This gives a closer 
analogy to the usual geometric notion of the antipodal map on $S^n$.

It is clear from the above construction that the triangulated $\RR P^4$ we constructed is invariant under the action 
of $S_6$ on $\{\ee_1,\ldots,\ee_6\}$. The action induced on the vertex set of this complex splits the vertices into 
two orbits, of size $6$ and $10$, corresponding to the quotients of $V_1\cup V_5$ and $V_3$ respectively. We refer to
this group of automorphisms as $\tilde{S}$, and it splits the $150$ facets of this complex into two orbits, of size 
$30$ and $120$. 

A quick comparison of the orbit representatives shows that the complex we have constructed is the same as the 
complex $\RR P^4_{16}$ in~\cite[p.77]{L}, available at \cite{Lw}. 

\subsection{A connection to the Witt design on $22$ points}

A convenient mnemonic to represent the facets of $\RR P^4$ in Construction~\ref{c1} is as follows. 

\begin{figure}
\begin{center}
 \begin{tikzpicture}[auto]
\node (A) at (120:4){A};
\node (B) at (60:4) {B} edge node[above,sloped]{$0123$}(A);
\node (C) at (360:4){C} edge node[below right,sloped]{$0456$}(A)
			edge node[above,sloped]{$0789$}(B);
\node (D) at (300:4){D} edge node[pos=.6, above,sloped]{$1489$}(A)
			edge node[above left,sloped]{$1567$}(B)
			edge node[below,sloped]{$2347$}(C);
\node (E) at (240:4){E} edge node[below right,sloped]{$2579$}(A)
			edge node[pos=.4, above ,sloped]{$2468$}(B)
			edge node[above left,sloped]{$1358$}(C)
			edge node[below,sloped]{$0369$}(D);
\node (F) at (180:4){F}   edge node[above,sloped]{$3678$}(A)
			  edge node[below right,sloped]{$3459$}(B)
			  edge node[pos=.6, above ,sloped]{$1269$}(C)
			  edge node[above left,sloped]{$0258$}(D)
			  edge node[below,sloped]{$0147$}(E);
      \end{tikzpicture}
 \caption{$K_6$ mnemonic for $16$-vertex $\RR P^4$}
\label{k6}
\end{center}
\end{figure}
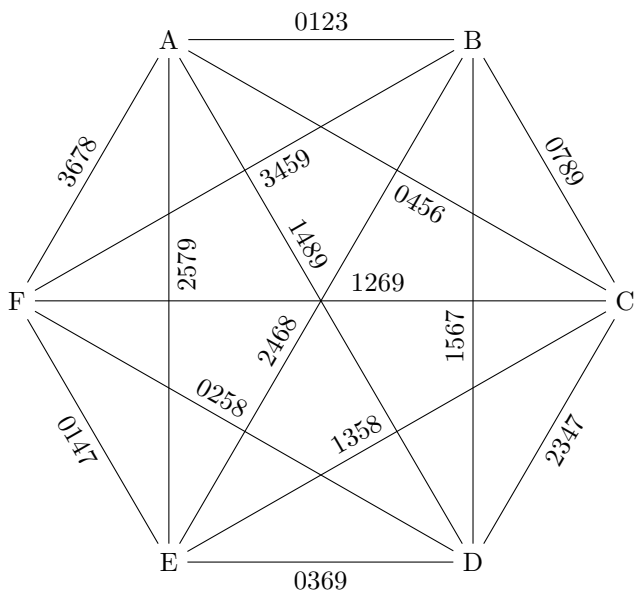

Label the vertices of the complete graph $K_6$ with the points of the $6$-vertex orbit under $\tilde{S}$. Now the remaining ten 
vertices correspond to the ten pairs of disjoint triangles, or \emph{bisections}, in $K_6$. We use these points to label 
the edges of $K_6$ as follows. Each edge is contained in four triangles, each of which is in turn contained in exactly 
one bisection. Give each edge a label consisting of the four bisections it is contained in. Henceforth, when we refer to  
$K_6$, we mean the $K_6$ labelled thus. We denote the set of elements of the $6$-vertex orbit by $\{A,\ldots,F\}$, and that 
of the $10$-vertex orbit by $\{0,\ldots,9\}$. See Figure~\ref{k6}.  

Now the small $\tilde{S}$-orbit of the set of facets can be read off Figure~\ref{k6} as a vertex-edge pair $(v,e)$ of $K_6$, where 
$v\in \{A,\ldots,F\}$ and $e$ is an edge-label of size $4$. The $\tilde{S}$-orbit of size $120$ is given by a triple 
$(v,v',e(v,v'')\setminus e(v,v'))$, where $v,v',v''\in \{A,\ldots,F\}$, and $e(u,v)$ denotes the $4$-label of the edge 
between $u$ and $v$. The correspondence between the facets of $S^4_{32}$ listed above and the facets read off $K_6$ above 
is straightforward.

It is fruitful to consider the properties of $K_6$ and its labels, from the point of view of combinatorial designs. Recall that 
a $t-(v,k,\lambda)$ \emph{design} $\mathbb{D}$ is a pair $(\mathbb{V},\mathbb{B})$, where $\mathbb{V}$ is a of size $v$, whose 
elements are called \emph{points}, and $\mathbb{B}$ is a set of $k$-subsets of $\mathbb{V}$ called \emph{blocks}, 
such that any $t$-subset of $\mathbb{V}$ is contained in exactly $\lambda$ elements of $\mathbb{B}$.

Let $\mathcal{V}$ (for vertices) denote the set $\{A,\ldots,F\}$, and $\mathcal{B}$ (for bisections) denote the set 
$\{0,\ldots,9\}$. Let $\mathcal{E}$ denote the set of fifteen edge-labels of $K_6$.

Our first observation is that $(\mathcal{B},\mathcal{E})$ a quasi-symmetric $2-(10,4,2)$ design. To see that this is a
$2$-design with $\lambda=2$, note that any pair of distinct bisections intersect in exactly two edges. To see that this 
design is quasi-symmetric, or that distinct blocks have two possible intersection sizes, namely $1$ and $2$, first note 
that the vertex-set of a pair of adjacent edges in $K_6$ forms a triangle, so is contained in a unique bisection, so their
labels intersect in one element, the bisection containing the triangle formed by their edges. Now if two edges are not 
adjacent, they are contained in exactly two bisections, and their labels intersect in those two bisections.\footnote{
This also shows that $\tilde{S}$ is the full automorphism group of $\RR P^4_{16}$. Counting the number of facets 
containing each vertex in $\mathcal{V}$ and $\mathcal{B}$, we see that the full automorphism group $G$ preserves these orbits.
But if $G>\tilde{S}\simeq S_6$, then $G$ can not act faithfully on $\mathcal{V}$, so there exists $g\in G\setminus \tilde{S}$ which 
fixes each block of $\mathcal{E}$. Then $g$ also fixes the intersections of these blocks, which implies $g$ fixes 
$\mathcal{B}$ pointwise. A contradiction.}
Also note that this design naturally extends to a $2-(16,6,2)$ design, with point-set $\mathcal{V}\cup\mathcal{B}$, and 
block-set $\widetilde{\mathcal{E}}:=\{\{u,v\}\cup e(u,v)|u,v\in\mathcal{V}\}\cup\{\mathcal{V}\}$. This design is also 
symmetric, as the number of blocks is the same as the number of points, and any two blocks intersect in $\lambda=2$ points. 

Further, recall that a \emph{perfect matching} or $1$-factor of a graph is a partition (if it exists) of its vertex set, 
into edges of the graph, and that being a complete graph on an even number of vertices, $K_6$ has perfect matchings, which 
we henceforth simply call matchings. Every edge of $K_6$ is contained  in exactly three matchings. The number of matchings 
in $K_6$ is $\frac{6!}{2!2!2!3!}=15$. The three pairwise disjoint edges in any matching have labels which intersect pairwise
in two elements of $\mathcal{B}$ each. But the bisections containing one pair of disjoint edges do not contain the third of 
these edges. So the union of three elements of $\mathcal{E}$ in a given matching is a subset of $\mathcal{B}$ of size $6$. 
So we can label each matching by the four elements of $\mathcal{B}$ not in the labels of any of its three edges. The set 
$\mathcal{M}$ (for matching) of these $4$-labels is the block set of another quasi-symmetric $2-(10,4,2)$ design, whose 
blocks intersect in one point if the corresponding matchings are disjoint and in two points if the matchings intersect 
in an edge. To see that $(\mathcal{B},\mathcal{M})$ is a $2$-design, note that deleting a pair of bisections from $K_6$ 
leaves the disjoint union of an edge and a $4$-cycle, which contains exactly two matchings. Also, to see that the intersection 
sizes are $2$ and $1$, note that the union of a pair of intersecting matchings is the disjoint union of an edge and a $4$-cycle, 
the complement of a pair of bisections, and that the union of two disjoint matchings is a hexagon, whose complement in $K_6$ 
contains exactly one bisection.

We introduce one final set of objects, the $1$-factorizations of $K_6$. Recall that a $1$-\emph{factorization} of a graph is 
a partition of its edge-set, where each block in the partition is a matching of the graph. The graph $K_6$ has six 
$1$-factorizations, which can be thought of as $5$-edge-colourings where the matchings are the colour-classes. 
We label these from the set $\{U,V,\ldots,Z\}=\mathcal{F}$ (for factorization).

It can also be seen that any matching is contained in exactly two $1$-factorizations, and that any two disjoint matchings 
determine a unique $1$-factorization. Now since any two $1$-factorizations can have at most one matching in common, and there are 
fifteen pairs of $1$-factorizations, any two $1$-factorizations, say $f,g$ intersect in a unique matching $m(f,g)$. This gives an 
extension of the design $(\mathcal{B},\mathcal{M})$ to a symmetric $2-(16,6,2)$ design $(\mathcal{B}\cup\mathcal{F},
\widetilde{\mathcal{M}})$, where  $\widetilde{\mathcal{M}}=\{\{f,g\}\cup m(f,g)|f,g\in\mathcal{F}\}\cup\{\mathcal{F}\}$. 

\begin{figure}
\begin{center}
 \begin{tikzpicture}
\node (A) at (120:2.8){A};
\node (B) at (60:2.8) {B} edge node[above,sloped]{\scriptsize $0123$}(A);
\node (C) at (360:2.8){C} edge node[below right,sloped]{\scriptsize $0456$}(A)
			edge node[above,sloped]{\scriptsize $0789$}(B);
\node (D) at (300:2.8){D} edge node[pos=.6, above,sloped]{\scriptsize $1489$}(A)
			edge node[above left,sloped]{\scriptsize $1567$}(B)
			edge node[below,sloped]{\scriptsize $2347$}(C);
\node (E) at (240:2.8){E} edge node[below right,sloped]{\scriptsize $2579$}(A)
			edge node[pos=.4, above ,sloped]{\scriptsize $2468$}(B)
			edge node[above left,sloped]{\scriptsize $1358$}(C)
			edge node[below,sloped]{\scriptsize $0369$}(D);
\node (F) at (180:2.8){F}   edge node[above,sloped]{\scriptsize $3678$}(A)
			  edge node[below right,sloped]{\scriptsize $3459$}(B)
			  edge node[pos=.6, above ,sloped]{\scriptsize $1269$}(C)
			  edge node[above left,sloped]{\scriptsize $0258$}(D)
			  edge node[below,sloped]{\scriptsize $0147$}(E);

\node [xshift=6cm](U) at (120:2.8){U};
\node [xshift=6cm](V) at (60:2.8) {V} edge node[above,sloped]{\scriptsize $5689$}(U);
\node [xshift=6cm](W) at (360:2.8){W} edge node[below right,sloped]{\scriptsize $1379$}(U)
			edge node[above,sloped]{\scriptsize $0249$}(V);
\node [xshift=6cm](X) at (300:2.8){X} edge node[pos=.6, above,sloped]{\scriptsize $1245$}(U)
			edge node[above left,sloped]{\scriptsize $0357$}(V)
			edge node[below,sloped]{\scriptsize $0168$}(W);
\node [xshift=6cm](Y) at (240:2.8){Y} edge node[below right,sloped]{\scriptsize $0348$}(U)
			edge node[pos=.4, above ,sloped]{\scriptsize $1278$}(V)
			edge node[above left,sloped]{\scriptsize $2356$}(W)
			edge node[below,sloped]{\scriptsize $4679$}(X);
\node [xshift=6cm](Z) at (180:2.8){Z}   edge node[above,sloped]{\scriptsize $0267$}(U)
			  edge node[below right,sloped]{\scriptsize $1346$}(V)
			  edge node[pos=.6, above ,sloped]{\scriptsize $4578$}(W)
			  edge node[above left,sloped]{\scriptsize $2389$}(X)
			  edge node[below,sloped]{\scriptsize $0159$}(Y);
  
\end{tikzpicture}
 \caption{$K_6$ and $K_6^*$}
\label{k6duals}
\end{center}
\end{figure}
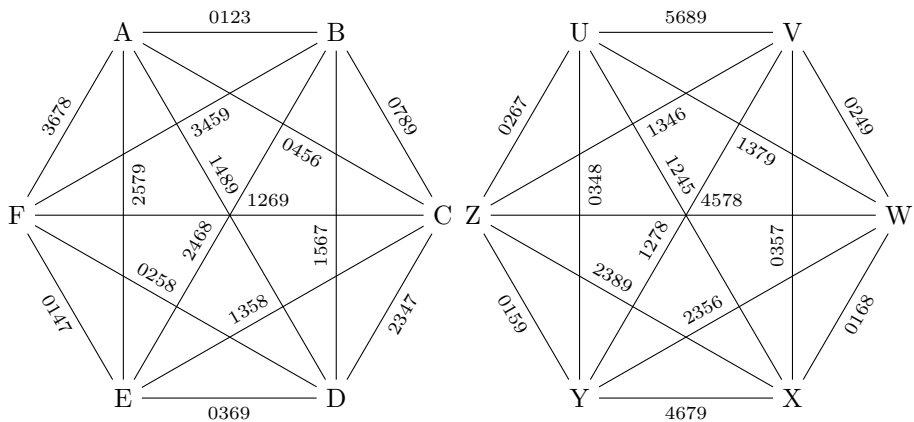

This can be represented by a ``dual $K_6$'', $K_6^*$ with vertices labelled by $\mathcal{F}$ and edges labelled by the elements of
$\mathcal{M}$, with $m(f,g)$ labelling the edge joining $f$ and $g$. Figure~\ref{k6duals} illustrates the two copies of 
$K_6$ with their edges labelled by the elements of $\mathcal{B}$. Note that the set $\mathcal{B}$ also indexes the bisections of 
$K_6^*$. This gives a correspondence between the $(3,3)-$partitions of $\mathcal{V}$ and those of $\mathcal{F}$.

We note the following type of triple-incidence between the sets $\mathcal{V},\mathcal{B},$ and $\mathcal{F}$.
Given any two bisections $b,b'$, each triangle in $b$ intersects one triangle in $b'$ in an edge $e$. 
This gives a partition of $\mathcal{V}$ into two edges and two points. The edges each correspond to 
the intersections of two triangles, one from $b$ and one from $b'$. The two points left over determine 
a third edge disjoint from the other two. Moreover, this edge, as a $4$-set in $\mathcal{E}$ contains 
neither $b$ nor $b'$. Also, since the pair $b,b'$ is contained in the labels of the other two edges, the
$4$-label in $\mathcal{M}$ of the $1-$factor $m$ composed of the three edges above is disjoint
from $\{b,b'\}$. Since $e\in m$, the $4$-labels of $e$ and $m$ are disjoint.

Similarly, given an incident edge-$1-$factor pair $(e,m)$, their $4$-labels are disjoint, and subtracting 
the union of these $4$-labels from $\mathcal{B}$ leaves two bisections such that the triangles of one 
intersects the triangles of the other in the two edges of $f\setminus e$. So we have a correspondence
between pairs of bisections and incident edge-$1-$factor pairs of $K_6$.

We are now ready to describe the $22$-point Witt design, the unique $3-(22,6,1)$ design first independently constructed 
by R.D. Carmichael and E. Witt in the 1930s~\cite{C,Wi}. To be accurate, we construct \emph{a} $3-(22,6,1)$ design, and 
make no claims as to its uniqueness. See, for example,~\cite{RS} for a proof of the uniqueness of $W_{22}$, which relies 
on the embeddability in it of a $16$-point biplane, which can be seen to be isomorphic to either $(\mathcal{V}\cup\mathcal{B},
\widetilde{\mathcal{E}})$ or $(\mathcal{B}\cup\mathcal{F},\widetilde{\mathcal{M}})$ described above. For our point set 
of $W_{22}$, we take $\mathcal{V}\cup\mathcal{B}\cup\mathcal{F}$. The blocks are of three types. The first two sets of 
blocks are $\widetilde{\mathcal{E}}$ and $\widetilde{\mathcal{M}}$, the block-sets of the two $2-(16,6,2)$ designs we 
described earlier. The third set of blocks is obtained as follows. For any incident edge-matching pair $(e,m)$ in $K_6$, 
let $v,v'\in\mathcal{V}$ be the vertices that make up $e$, and let $f,f'\in\mathcal{F}$ be the $1$-factorizations that 
intersect in $m$. Let $\underline{e}$ and $\underline{m}$ be the $4$-sets corresponding to $e$ and $m$ in $\mathcal{E}$ 
and $\mathcal{M}$ respectively. Now for each of the $45$ incident pairs $(e,m)$, we take the block $\{v,v',f,f'\}\cup 
(\mathcal{B} \setminus(\underline{e}\cup \underline{m}))$. Call the set of such blocks $\mathcal{EM}$.

\begin{thm}
The pair $(\mathcal{V}\cup\mathcal{B}\cup\mathcal{F},\widetilde{\mathcal{E}}\cup\widetilde{\mathcal{M}}
\cup\mathcal{EM})$ defined above is a $3-(22,6,1)$ design.
\end{thm}
\begin{proof}
There are $77$ blocks, each of size $6$. As there are $22$ points, and since $\binom{22}{3}=77\times\binom{6}{3}$, it is
enough to show that every $3$-subset of the point set appears in some block.

Consider the $3$-subsets of $\mathcal{V}\cup\mathcal{B}\cup\mathcal{F}$. As $\mathcal{V}$ and $\mathcal{F}$
are themselves blocks, any $3$-subset of either is in a  block. There are $\binom{10}{3}=120$ subsets of 
$\mathcal{B}$ of size $3$. Consider the $3$-subsets of $\mathcal{E}$ or $\mathcal{M}$. If we can show that 
no two sets of these forms have a $3$-subset in common, it will follow that there are $2\times 15\times 4=120$ 
such sets, and that every $3$-subset of $\mathcal{B}$ is in exactly one block. Two elements of $\mathcal{E}$ intersect 
in at most two points of $\mathcal{B}$. Similarly, two elements of $\mathcal{M}$ intersect in at most two points of 
$\mathcal{B}$. Now consider the intersection of an element $\underline{e}\in\mathcal{E}$ with $\underline{m}\in\mathcal{M}$.
If the corresponding edge and matching of $K_6$ are incident, they do not intersect at all. Otherwise, the edge 
corresponding to $\underline{e}$ has its vertices in two (disjoint) edges of $\underline{m}$, say $uu'$ and $vv'\subset 
\mathcal{V}$. (We write a set as a string from here on, for convenience and brevity.) Now since the labels of incident 
edges have one element in common, $|e(u,v)\cap e(u,u')| =|e(u,v)\cap e(v,v')|=1$. Also since $u'\neq v'$, $e(u,v)\cap 
e(u,u')\neq e(u,v)\cap e(v,v')$. Since the labels of the third edge in the matching corresponding to $\underline{m}$ are 
contained in $e(u,u')\cup e(v,v')$, we have $|\underline{e}\cap \underline{m}|=|\underline{e}\cap (\mathcal{B}\setminus
(e(u,u')\cup e(v,v')))|=2$. So all the $3$-subsets of the $4$-labels of edges and matchings are distinct, and there are 
$120$ such sets, so each $3$-subset of $\mathcal{B}$ is in exactly one block.

Now consider a $3$-set consisting of two elements of $\mathcal{V}$ and one element of $\mathcal{F}$, say $vv'f$. Since 
$f$ contains exactly one matching containing the edge $vv'$, this $3$-set is contained in exactly one block of $\mathcal{EM}$.
Similarly, given a $3$-set of the form $vff'$, the $1$-factorizations $f$ and $f'$ intersect in a unique matching, 
which is a partition of the vertex set of $K_6$. So $vff'$ is in exactly one element of $\mathcal{EM}$.

A $3$-set of the form $vv'b$ is of one of the following types. If $b\in e(v,v')$, then $vv'b$ is contained (again, 
in fact, exactly once) in a block of $\widetilde{\mathcal{E}}$. If $b\notin e(v,v')$, then $v,v'$ are in different 
triangles of $b$, and $v,v'$ together with the remaining  edges of these triangles forms a matching whose label does
not contain $b$. So $vv'b$ is contained in a block of $\mathcal{EM}$.

Similarly, if $b\in m(f,f')$, then $ff'b$, is contained in exactly one block of $\mathcal{M}$. If $b\notin m(f,f')$,
then it is an element of a $4$-label of two of the edges in the matching $m$ that $f$ and $f'$ intersect in. Let $e$ be the 
third edge  of this matching. Then $ff'b$ is in the block of $\mathcal{EM}$ corresponding to $(e,m)$.

Now if a $3$-set is of the form $vbb'$, we have two possibilities. Consider the two blocks of 
$\mathcal{E}$ containing $bb'$ . If $v$ is in either of the corresponding edges of $K_6$, then $vbb'$ is in 
$\widetilde{\mathcal{E}}$. The two vertices not in either of these edges, form the third edge of the matching 
containing the two edges whose labels contain $bb'$. So if $v$ is on this edge, $vbb'$ is in a block of 
$\mathcal{EM}$. Similarly a $3$-set of the form $fbb'$ is in a block of either $\widetilde{\mathcal{M}}$ or $\mathcal{EM}$. 

The only remaining type of $3$-set is of the form $vfb$. These can only be contained in blocks of $\mathcal{EM}$. 
Now there are $45$ blocks in $\mathcal{EM}$, each corresponding to a pair of elements in $\mathcal{B}$. So each 
element $b$ of $\mathcal{B}$ is in nine blocks of $\mathcal{EM}$, corresponding to the nine edges of $K_6$ not labelled 
by $b$. Let $v\in\mathcal{V}$ and let the bisection of $K_6$ corresponding to $b$ be $vv'v'',uu'u''$. 
Then $v\in\mathcal{V}$ is in three edges whose labels do not contain $b$, each of which determine exactly one matching 
whose label does not contain $b$. All three matchings intersect in the edge $v'v''$. Since a $1$-factorization is a partition 
of the edge set of $K_6$, no $1$-factorization contains more than one of these matchings. So the pairs of $1$-factors which 
intersect to give each of these matchings partition the vertex set of $K_6^*$. In other words, the two elements of 
$\mathcal{F}$ in each of the three blocks of $\mathcal{EM}$ containing a pair $vb$ are all distinct. So any triple $vfb$ is 
in a block of $\mathcal{EM}$. This completes the proof.
\end{proof}
The arguments above also illustrate the duality between $K_6$ and $K_6^*$. All properties of the vertices and edges of $K_6$ 
have analogues in its $1$-factorizations and matchings, or the vertices and edges of $K_6^*$. Applying the dual construction on 
$K_6^*$ simply recovers $K_6$. Let $S_6'$ be the group of permutations of $\mathcal{F}$ induced by the permutation group $S_6$ 
of $\mathcal{V}$ via permutations of $\mathcal{B}$. The actions of $S_6$ and $S_6'$ on $\mathcal{V}$ and $\mathcal{F}$ 
respectively are \emph{dual} or non-conjugate to each other, i.e., there is no one-to-one map between the elements of these 
two sets which send either group to the other. Abstractly, any isomorphism between these two sets corresponds to an outer 
automorphism of $S_6$. See~\cite[Chapter 6]{CvL} for a more thorough treatment and several applications. In particular, we note 
that our construction of $W_{22}$ has interesting parallels with the construction of the $5-(12,6,1)$ Witt design 
on~\cite[p.86]{CvL}.

The reader familiar with the Witt designs and their automorphism groups will recognise that our description of $W_{22}$ 
can be recovered from its usual forms by fixing any pair of its disjoint blocks and mapping it to $(\mathcal{V},\mathcal{F})$. 
The ten remaining points correspond exactly to $(3,3)$-partitions of $\mathcal{V}$, and also of $\mathcal{F}$. 

We take this occasion to note early, if not original, appearances of the various objects that feature in our construction in the 
literature. The remarkable duality between the vertices and edges of $K_6$ and its $1-$factors and $1-$factorizations was 
noted by J.J. Sylvester in a paper from 1844, republished in~\cite{S}. The origins of the (isomorphism class of) $16$-point 
biplanes are historic, and can be traced back to Kummer's $16_6$ configuration. The combinatorial aspects of this configuration
are studied in great depth in the classic treatise of R.W.H.T. Hudson,~\cite[Sections 5,9\footnote{We would like to draw special 
attention to this section, which demonstrates, for any choice of bisection $b$ of $K_6$, a correspondence between the set 
$\mathcal{V}\cup\mathcal{B}\setminus\{b\}$ consisting of its vertices and the remaining bisections, and the set of the edges of $K_6$. When 
the $1-$factors of $K_6$ are added to both sets, the latter set corresponds to the point-set of the $21$-point projective plane 
$\Pi_4$. The bisection $b$ now corresponds to the one point which is used to extend $\Pi_4$ in the Witt-L\"uneberg construction 
of $W_{22}$.},23,24,25,26]{H}. The observation that the actions of $S_6$ and $S_6'$ on $\mathcal{V}$ and $\mathcal{F}$ 
respectively induce the same group of permutations on the set $\mathcal{B}$ of bisections of $K_6$ and $K_6^*$ appears in~\cite{E}. 

We also point out some similarities and differences between our construction and two earlier ones. In~\cite[Table 8.2]{M}, D.M. 
Mesner, at the time unaware of the preexistence of $W_{22}$ in the group theoretic literature, lists the blocks of a ${3-(22,6,1)}$ 
design discovered by empirical search. In his proof of the uniqueness of this design~\cite[Theorem 8.7]{M}, he notes and relies on 
the facts that every block is disjoint from $16$ other blocks, that every pair of points is contained in five blocks, and that the 
possible intersection sizes of blocks are $0$ or $2$. These observations enable him to fix an initial block, list up to 
relabellings the $60$ blocks that intersect the initial blocks, then show that the remaining $16$ blocks disjoint from the initial 
block are also fixed, noting that these are the blocks of a $2-$design. See~\cite{KW} for a historical account and survey of these 
and further implications of Mesner's work. 

In~\cite{RS}, N.N. Roghelia and S.S. Sane, independently of Mesner, but familiar with the work of Witt as expounded by 
H. L\"uneberg~\cite{Lu}, prove the existence and uniqueness of $W_{22}$ based on the classification of $16$-point biplanes by the 
number of their \emph{ovals}, or sets of four points, no three of which are in a block. They show that the unique $16$-point biplane 
with $60$ ovals is uniquely embeddable in $W_{22}$, by way of the $2-(16,4,3)$ design with these ovals as blocks, and six new points.
We also note that though Roghelia and Sane's construction refer to the complete graph $K_6$, this graph is indexed by blocks, unlike 
our $K_6$.

We would like to think of the construction of Mesner as starting with $6+16$ points and constructing $1+60+16$ blocks on them. 
Roghelia and Sane proceed in the opposite direction, starting with $16+6$ points and constructing $16+60+1$ blocks. Within this 
point of view, our construction starts with $6+10+6$ points and constructs $16+45+16$ blocks on these. At the time of discovering 
our construction, this author was only aware of the better known construction of $W_{22}$ via the extension of the projective 
plane $\Pi_4$, and as embedded in the larger design $W_{24}$. Despite various aspects of our construction already existing in 
the literature, we find our view of $W_{22}$ where two biplanes coexist, mediated by the bisections of the dual pair of $K_6$ and 
$K_6^*$, to be of some interest.

We also remark that the connection between the $16$-vertex $\RR P^4$ we constructed earlier and the above $3$-design is not 
simply restricted to the $2-(16,6,2)$ design $(\mathcal{V}\cup\mathcal{B},\widetilde{\mathcal E})$.The $\tilde{S}$-orbit of 
the facets of our $\RR P^4$ split naturally into fifteen sets of size eight, each containing a pair of vertices $v,v'$ in 
$\mathcal{V}$. The  closure of any of these $8$-sets as facets, is the join of the edge $[v,v']$ with its link in our 
$\RR P^4$, an octahedron with vertices from $\mathcal{B}$. Any octahedron is determined by three pairs of opposite vertices
on each ``axis'', say $b_1b_1',b_2b_2',b_3b_3'$ in this case. Then each of the $15\times 3=45$ quadruples 
$vv'b_ib_i', 1\leq i\leq 3$ is contained in a block of $\mathcal{EM}$.
\footnote{These are also ovals of the biplane $(\mathcal{V}\cup\mathcal{B},\widetilde{\mathcal{E}})$ used in~\cite{RS}. 
The remaining $15$ ovals are the blocks of the design $(\mathcal{B},M)$.} 
Also if we take the blocks of $W_{22}$ and delete the elements of one block from all the others, the 
remaining blocks split into a set of sixteen blocks of size $6$, and sixty blocks of size $4$. The set $\mathbb{B}_6$ with 
sixteen blocks of size $6$ form the block-set of a symmetric $2-(16,6,2)$ design, which corresponds to the block-set 
$\widetilde{\mathcal{E}}$ we started with. Now if we fix a block $B$ of $\mathbb{B}_6$, this further splits the set 
$\mathbb{B}_4$ of $4$-sets into $15$ sets that are disjoint from it, (i.e., $\mathcal{M}$), and $45$ sets that intersect $B$ 
in two points, in correspondence with $\mathcal{EM}$. So we can think of each block of $W_{22}$ as sitting in the centre of a 
configuration of $16$ copies of $\RR P^4_{16}$.

\section{Further Constructions}\label{sec:fc}
It would be of interest to know if there is a general ``algorithm'' to construct minimal, or even smaller-than-known 
triangulations of real projective spaces. Our first construction of triangulated $\RR P^4$, though short and 
straightforward, has some exceptional properties which may well be the results of numerical coincidences, and 
do not offer much hope of analogous constructions in higher dimensions. We provide two more ways of constructing a 
$16$-vertex $\RR P^4$ which are easier to generalise.

\begin{rk}
It must be borne in mind here that even though our choice of notation in the following constructions represents the
vertices of $n$-dimensional complexes as points in $\RR^n$, the objects we construct are purely abstract simplicial 
complexes, which we do not need to view as embedded in $\RR^N$ for any $N$. Indeed, they most definitely do not embed
in $\RR^n$. Our choice of notation is motivated by ease of handling and conceptual visualization.
\end{rk}
 
\subsection{Constructions using cross-polytopes and hypercubes} 
We explore the possibility of constructing triangulated real projective $n$-space, $\RR P^n$ in the following way. Take an
$n$-dimensional cross-polytope and triangulate its interior, possibly by adding an extra point $\mathbf{0}$. Say we denote 
the vertices of the cross-polytope $C^n$ by the vectors $\pm\ee_i\in\RR^n$, where $\ee_i$ is the $i^{\mathrm{th}}$ elementary
vector, $1\leq i\leq n$. Then we add $2^n$  simplices of the form 
$[\varepsilon_1\ee_1,\varepsilon_2\ee_2,\ldots,\varepsilon_n\ee_n, \sum_{i=1}^n\varepsilon_i\ee_i]$, where 
$\varepsilon=(\varepsilon_1,\varepsilon_2,\ldots,\varepsilon_n)\in\{\pm1\}^n$. Call the new vertex set $V$. We have 
$V=C\sqcup Q$, where $C$ is the vertex set of the triangulated $C^n$ and $Q=\{q_{\varepsilon}|\varepsilon\in\{\pm1\}^n\}=
\{\sum_{i=1}^n \varepsilon_i\ee_i|\varepsilon_i=\pm1,1\leq i\leq n\}$.

Now we consider subsimplices of $\partial C^n$ going down in dimension, and triangulate the links of each without adding
any more vertices. Our goal is to end up with the boundary of the link of a vertex in $C$ as a triangulated $2^{n-1}$-vertex
$n-2$-sphere. We do this subject to the following conditions. First, for every facet
$$[\varepsilon_{i_1}\ee_{i_1},\varepsilon_{i_2}\ee_{i_2},\ldots,\varepsilon_{i_k } \ee_{i_k},q_{ \varepsilon^{j_1}},
q_{\varepsilon^{j_2}},\ldots,q_{\varepsilon^{j_{n-k+1}}}]$$
 in the complex, its ``opposite'' facet 
$$[-\varepsilon_{i_1}\ee_{i_1},-\varepsilon_{i_2}\ee_{i_2},\ldots,-\varepsilon_{i_k }\ee_{i_k},q_{-\varepsilon^{j_1}},
q_{-\varepsilon^{j_2} } ,\ldots,q_{-\varepsilon^{j_{n-k+1}}}]$$
is also in the complex. Second, no vertex $q_{\varepsilon}$ in $Q$ is joined to its ``opposite'' vertex 
$q_{-\varepsilon}=-q_{\varepsilon}$. In other words $[q_{\varepsilon},q_{-\varepsilon}]$ is not an edge of the complex.
Third, if a vertex $u\in V$ is joined to another vertex $v\in V$, i.e, if $[u,v]$ is an edge of the complex, then $[u,-v]$
is not an edge of the complex. These conditions, equivalent to the existence of a link-separating involution on the complex, 
allow us to apply the identification map $q_{\varepsilon}\sim q_{-\varepsilon}$ on $Q$, leaving us with $2n$ suspended 
$S^{n-2}$. We then triangulate the interior of these $2^{n-1}+2$-vertex $n-1$-spheres, to get a triangulation of $\RR P^n$.
\begin{ex}
It is easily seen that the paradigm outlined above can be used to construct $\RR P^2_6$ 
as follows. See Figure ~\ref{rp2}. We triangulate the square with vertices $\pm\ee_1,\pm\ee_2$
by joining $+\ee_1$ with $-\ee_1$, into triangles $[\ee_1,-\ee_1,\ee_2],[\ee_1,-\ee_1,-\ee_2]$.

Next we add the triangles
$\pm[\ee_1,\ee_2,\ee_1+\ee_2],\pm[\ee_1,-\ee_2,\ee_1-\ee_2]$. The links of the
vertices $\pm\ee_1$ now have boundaries $\{\pm[\ee_1+\ee_2],\pm[\ee_1-\ee_2]\}$
respectively and the boundaries of the links of $\pm\ee_2$ have boundaries
$\{\pm[\ee_1+\ee_2],\pm[-\ee_1+\ee_2]\}$. Now we apply the map.
$q_{\varepsilon}\sim q_{-\varepsilon}$.
We triangulate the $1$-sphere containing $\pm\ee_2$ by adding the triangles
$[+\ee_2,-\ee_2,\overline{\ee_1+\ee_2}],[+\ee_2,-\ee_2,\overline{\ee_1+\ee_2}]$.
Since the link of $[+\ee_1,-\ee_1]$ is already a $0$ sphere in our complex, we
triangulate the remaining square as
$[+\ee_1,\overline{\ee_1+\ee_2},\overline{\ee_1-\ee_2}],[-\ee_1,\overline{
\ee_1+\ee_2},\overline{\ee_1-\ee_2}]$. This gives us $\RR P^2_6$. 

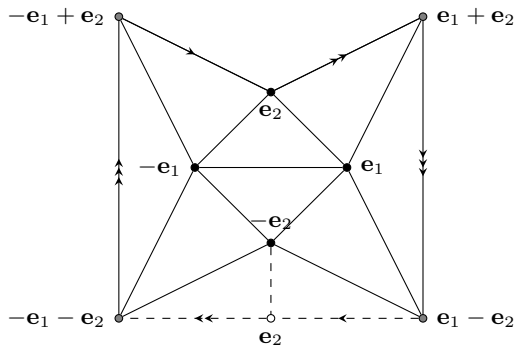
\begin{figure}%[b]
\begin{center}
\begin{tikzpicture}[>=stealth,vx/.style={circle,inner sep=0pt,minimum size=1mm,draw}]
\node (1) 	at 	(1,0) [vx] [fill,label=right:\small $\ee_1$]{}; 
\node (-1) 	at 	(-1,0) [vx] [fill,label=left:\small$-\ee_1$]{} ; 
\node (2)	at 	(0,1) [vx] [fill,label=below:\small$\ee_2$]{} ; 
\node (-2)	at 	(0,-1) [vx] [fill,label=above:\small$-\ee_2$]{} ; 
\node (++) 	at 	(2,2) [vx] [fill=black!50,label=right:\small$\ee_1+\ee_2$]{}; 
\node (-+)	at 	(-2,2) [vx] [fill=black!50,label=left:\small$-\ee_1+\ee_2$]{}; 
\node (+-)	at 	(2,-2) [vx] [fill=black!50,label=right:\small$\ee_1-\ee_2$]{}; 
\node (--)	at 	(-2,-2) [vx] [fill=black!50,label=left:\small$-\ee_1-\ee_2$]{}; 
\node (down) at (0,-2) [vx] [label=below:\small$\ee_2$]{} ;
\draw (1)--(-1)--(2)--(1)--(-2)--(-1);
\draw (+-)--(1)--(++)--(2)--(-+)--(-1)--(--)--(-2)--(+-);
\draw [decoration={markings, mark=at position .5 with {\arrow{>>}},mark=at position .53 with {\arrow{>}}},
		  postaction={decorate}](++)--(+-);
\draw [decoration={markings, mark=at position .5 with {\arrow{>>}},mark=at position .53 with {\arrow{>}}}, 
		  postaction={decorate}](--)--(-+);
\draw [decoration={markings, mark=at position .75 with {\arrow{>>}}, mark=at position .25 with {\arrow{>}}},
		  postaction={decorate},dashed]  (+-)--(down)--(--);
\draw [decoration={markings, mark=at position .25 with {\arrow{>}}, mark=at position .75 with {\arrow{>>}}},
		  postaction={decorate}] (-+)--(2)--(++);
\draw [dashed]	 (-2)--(down);
\end{tikzpicture}
\caption{$\RR P^2$ triangulated with squares}
\label{rp2}
\end{center}
\end{figure}
\end{ex}

\begin{ex}
In order to construct $\RR P^3$ by the same approach, start with the octahedron $C^3$ spanned by the 
points $\pm\ee_i$, where $i=1,2,3$. We can triangulate the interior of the octahedron by taking the cone
over its boundary at the point $\mathbf{0}$. This gives us eight tetrahedra of the form 
$$[\mathbf{0},\varepsilon_1\ee_1,\varepsilon_2\ee_2,\varepsilon_3\ee_3].$$

The boundary of the octahedron consists of the eight triangles of the form 
$[\varepsilon_1\ee_1,\varepsilon_2\ee_2,\varepsilon_3\ee_3]$, where $\varepsilon_i=\pm1$ for $i=1,2,3$. Now add a set 
$Q$ of eight new ``outer'' vertices ${q_{\varepsilon}=\varepsilon_1\ee_1+\varepsilon_2\ee_2+\varepsilon_3\ee_3}$ for
each $\varepsilon=(\varepsilon_1,\varepsilon_2,\varepsilon_3)\in\{\pm1\}^3$, by taking the eight tetrahedra of the form
$$[\varepsilon_1\ee_1,\varepsilon_2\ee_2,\varepsilon_3\ee_3,\varepsilon_1\ee_1+\varepsilon_2\ee_2+\varepsilon_3\ee_3].$$ 
The boundary of this complex is a triangulated $S^2$ with $f$-vector $[14,36,24]$. 

Now consider the link of an edge of $C^3$. The link of $[\varepsilon_i\ee_i,\varepsilon_j\ee_j]$ is the path 
$[\varepsilon_i\ee_i+\varepsilon_j\ee_j-\ee_k,-\ee_k],[-\ee_k,\mathbf{0}],[\mathbf{0},\ee_k],
[\ee_k,\varepsilon_i\ee_i+\varepsilon_j\ee_j+\ee_k]$, where $\{i,j,k\}=\{1,2,3\}$. Its boundary consists of 
the two points $\varepsilon_i\ee_i+\varepsilon_j\ee_j\pm\ee_k$. Close the boundary of
$[\varepsilon_i\ee_i,\varepsilon_j\ee_j]$ by adding the tetrahedron
$$[\varepsilon_i\ee_i,\varepsilon_j\ee_j,{\varepsilon_i\ee_i+\varepsilon_j\ee_j-\ee_k,
\varepsilon_i\ee_i+\varepsilon_j\ee _j+\ee_k}].$$ 
This gives twelve new tetrahedra, and the boundary of the new complex is still a triangulated $S^2$ with fourteen 
vertices. But the boundary of the link of a vertex $\varepsilon_i\ee_i$ of $C^3$ is now the boundary of a square, 
whose vertices are all its neighbours in $Q$. See Figure ~\ref{rp3} for an illustration of the neighbourhood of 
$\varepsilon_i\ee_i$. Also note that the subcomplex spanned by $Q$ is the set of edges of the $3$-dimensional cube $Q^3$.

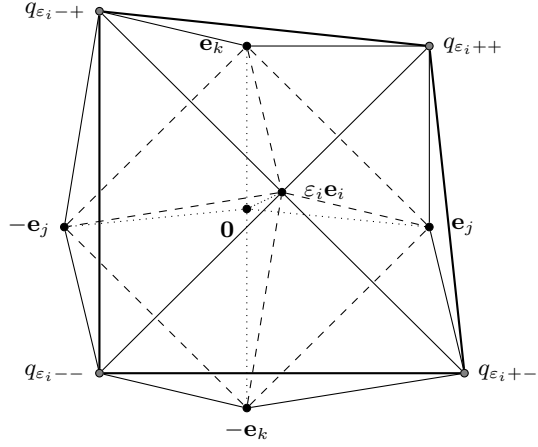
\begin{figure}%[b]
\begin{center}
\begin{tikzpicture}[scale=1.2,vx/.style={circle,inner sep=0pt,minimum size=1mm,draw}]
\node (0) 	at 	(0,0.2,1) [vx] [fill,label=below left:\small $\mathbf{0}$]{}; 
\node (1) 	at 	(2,0,1) [vx] [fill,label=right:\small $\ \ee_j$]{}; 
\node (-1) 	at 	(-2,0,1) [vx] [fill,label=left:\small$-\ee_j$]{} ; 
\node (2)	at 	(0,2,1) [vx] [fill,label=left:\small${\ee_k\ }$]{} ; 
\node (-2)	at 	(0,-2,1) [vx] [fill,label=below:\small$-\ee_k$]{} ; 
\node (3)	at 	(0,0,0) [vx] [fill,label=right:\small $\ \varepsilon_i\ee_i$]{} ; 
\node (+++) 	at 	(2,2,1) [vx] [fill=black!50,label=right:\small $q_{\varepsilon_i++}$]{}; 
\node (-++)	at 	(-2,2) [vx] [fill=black!50,label=left:\small $q_{\varepsilon_i-+}$]{}; 
\node (+-+)	at 	(2,-2) [vx] [fill=black!50,label=right:\small $q_{\varepsilon_i+-}$]{}; 
\node (--+)	at 	(-2,-2) [vx] [fill=black!50,label=left:\small $q_{\varepsilon_i--}$]{};
\draw[dotted]		(1)--(0)--(-1) (2)--(0)--(-2);
\draw[densely dotted] (0)--(3);
\draw[dashed]		(-2)--(3)--(1)--(2)--(3)--(-1)--(-2)--(1) (-1)--(2);
\draw			(1)--(+++)--(2)--(-++)--(-1)--(--+)--(-2)--(+-+)--(1) (+++)--(3)--(+-+) (-++)--(3)--(--+) ;
\draw[thick]			(+++)--(+-+)--(--+)--(-++)--(+++);
\end{tikzpicture}
\caption{Neighbourhood of $\varepsilon_i\ee_i$ in a $3$-ball before quotienting.}
\label{rp3}
\end{center}
\end{figure}

We can now identify $\sum_{i=1}^3\varepsilon_i\ee_i$ with the point $-\sum_{i=1}^3\varepsilon_i\ee_i$, and close the 
boundary of the link of $\varepsilon_i\ee_i$ by taking its cone at the point $-\varepsilon_i\ee_i$. That is, for each
pair $\pm\ee_i$, we take the four tetrahedra 
$$[\varepsilon_i\ee_i,-\varepsilon_i\ee_i,\pm\overline{\ee_j+\ee_k+\varepsilon_i\ee _i},
\overline{\ee_j+\ee_k+\varepsilon_i\ee_i}]$$
The link of the vertex $\varepsilon_i\ee_i$ is now the triangulated $8$-vertex
$S^2$ with facets
\begin{multline*}
[\mathbf{0},\pm\ee_j,\pm\ee_k],\pm[\ee_j,\ee_k,\overline{
\ee_j+\ee_k+\varepsilon_i\ee_i}],\pm[\ee_j,-\ee_k,\overline{
\ee_j-\ee_k+\varepsilon_i\ee_i}],\\
[\pm\ee_j,\overline{\ee_j+\ee_k+\varepsilon_i\ee_i},\overline{
\ee_j-\ee_k+\varepsilon_i\ee_i }],[\overline{ \ee_j+\ee_k+\varepsilon_i\ee_i
} ,\overline{\ee_j-\ee_k+\varepsilon_i\ee_i}, -\varepsilon_i\ee_i]
 \end{multline*}
The above complex is an $11$-vertex triangulation of $\RR P^3$. This complex is the same as the minimal
$\RR P^3_{11}$ described by Walkup in~\cite{W}, as the antipodal quotient of a $22$-vertex $S^3$.
\end{ex}

We now tackle $\RR P^4$.

\begin{cons} \label{c2}
We start with a $4$-dimensional (solid) hyperoctahedron $C^4$, given by the convex hull of $\{\pm \ee_1, \pm \ee_2, 
\pm \ee_3,\pm \ee_4 \}$. We triangulate $C^4$ by joining the vertices $+\ee_1$ and $-\ee_1$. The resulting complex is
a set of eight $4$-simplices which can be visualized as the join of the line segment $[-\ee_1,+\ee_1]$ with the 
boundary of the octahedron spanned by $\{\pm \ee_2, \pm \ee_3,\pm \ee_4 \}$.

The boundary of this triangulated $C^4$ is just the boundary $\partial C^4$ of $C^4$, which is a triangulated 
$3$-sphere with $f$-vector $[8, 24, 32, 16]$. Now we take the cone over each of the $16$ facets of this boundary 
with a different point. That is, for each facet $[\varepsilon_1\ee_1,\varepsilon_2\ee_2,\varepsilon_3\ee_3,
\varepsilon_4\ee_4]$ of $\partial C^4$, take the cone over this facet at the point 
$q_{\varepsilon}=\sum_{i=1}^4\varepsilon_i\ee_i$, where 
$\varepsilon=(\varepsilon_1,\varepsilon_2, \varepsilon_3,\varepsilon_4)\in\{\pm1\}^4$. 
This gives sixteen such $4$-simplices, and the boundary now has $24$ vertices, $16\times 4 + 24=88$ edges, 
$16\times \binom{4}{2}+32=128$ triangles and $16\times 4=64$ tetrahedra. Denote this triangulation of $S^3$ by $X^{(1)}$.

Now consider the link of each triangle $[\varepsilon_i\ee_i,\varepsilon_j\ee_j,\varepsilon_k\ee_k]$
of $\partial C^4$ in $X^{(1)}$. These are of two kinds. If $1\notin\{i,j,k\}$, then the link of  
$[\varepsilon_i\ee_i,\varepsilon_j\ee_j,\varepsilon_k\ee_k]$ is 
$$[-\ee_1+\varepsilon_i\ee_i+\varepsilon_j\ee_j+\varepsilon_k\ee_k,-\ee_1],
[-\ee_1,+\ee_1],[+\ee_1,+\ee_1+\varepsilon_i\ee_i+\varepsilon_j\ee_j+
\varepsilon_k\ee_k].$$
Now suppose $1\in\{i,j,k\}$, then the link of the triangle is
$${[-\ee_l+\sum_{i,j,k}\varepsilon_\alpha\ee_\alpha , -\ee_l ] , [ -\ee_l,-\varepsilon_1\ee_1],
[-\varepsilon_1\ee_1,+\ee_l] , [+\ee_l,+\ee_l+\sum_{i,j,k}\varepsilon_\alpha\ee_\alpha],}$$
where $l$ is the coordinate in $\{1,2,3,4\}\setminus\{i,j,k\}$.
In either case the endpoints of the link of the triangle are the two points $q_{\varepsilon^{+l}}$ and 
$q_{\varepsilon^{-l}}$ corresponding to the two tetrahedra containing it in $C^4$. 

We can now close the links of the triangle $[\varepsilon_i\ee_i,\varepsilon_j\ee_j,\varepsilon_k\ee_k]$ by adding
the $4$-simplices
$$[\varepsilon_i\ee_i,\varepsilon_j\ee_j,\varepsilon_k\ee_k, \sum_{i,j,k}\varepsilon_{\alpha}\ee_{\alpha}-\ee_l
,\sum_{i,j,k} \varepsilon_{\alpha}\ee_{\alpha}+\ee_l]$$
where $\{i,j,k,l\}=\{1,2,3,4\}$. We have added $32$ such $4$-simplices, and the boundary of the new 
complex is a triangulated $S^3$ with $24$ vertices, $88+32=120$ edges, $128-32+(32\times 3)=192$ triangles, and $64+32=96$
tetrahedra. Call the boundary $X^{(2)}$. 

Observe now that the subcomplex spanned by the subset of vertices 
$Q=\{q_{\varepsilon}|\varepsilon\in\{\pm 1\}^4\}$ is the $1$-skeleton of a $4$-dimensional hypercube, 
which is the dual of $C^4$. 

Now consider the links of the edges of $C^4$. The link of the edge $[\varepsilon_i\ee_i,\varepsilon_j\ee_j]$, when 
$1\notin\{i,j\}$ has eight vertices, namely $\pm\ee_k,\pm\ee_1, \varepsilon_i\ee_i+\varepsilon_j\ee_j\pm\ee_k\pm\ee_1$, 
where $\{i,j,k\}=\{2,3,4\}$. These vertices can be seen as forming the corners of a $4$-sided antiprism whose opposite 
squares are $(+\ee_1,+\ee_k,-\ee_1,-\ee_k)$ and 
$(\varepsilon_i\ee_i+\varepsilon_j\ee_j+\ee_k+\ee_1, \varepsilon_i\ee_i+\varepsilon_j\ee_j+\ee_k-\ee_1, 
\varepsilon_i\ee_i+\varepsilon_j\ee_j-\ee_k-\ee_1, \varepsilon_i\ee_i+\varepsilon_j\ee_j-\ee_k+\ee_1)$. 
All edges and triangles of this antiprism are faces of $X^{(2)}$. The (cyclically ordered) square 
$(+\ee_1,+\ee_k,-\ee_1,-\ee_k)$ is triangulated by the edge $[\ee_1,-\ee_1]$, and no triangle with vertices 
from the the other square is a face of $X^{(2)}$. The link of the edge $[\varepsilon_1\ee_1,\varepsilon_j\ee_j]$
is almost the same, with the only difference being that the square $(+\ee_k,+\ee_l,-\ee_k,-\ee_l)$ is triangulated 
by taking the cone of its boundary at $-\varepsilon_1\ee_1$. See Figure~\ref{antiprism}.
\begin{figure}
\begin{center}
\begin{tikzpicture}[vx/.style={circle,inner sep=0pt,minimum size=1mm,draw}]
\node (1) 	at 	(1,0) [vx] [fill,label=right:\small $\ee_l$]{}; 
\node (-1) 	at 	(-1,0) [vx] [fill,label=left:\small$-\ee_l$]{} ; 
\node (2)	at 	(0,1) [vx] [fill,label=above:\small$\ee_k$]{} ; 
\node (-2)	at 	(0,-1) [vx] [fill,label=below:\small$-\ee_k$]{} ; 
\node (++) 	at 	(2,2) [vx] [fill=black!50,label=above:\small$\varepsilon_i\ee_i+\varepsilon_j\ee_j+\ee_k+\ee_l$]{}; 
\node (-+)	at 	(-2,2) [vx] [fill=black!50,label=above:\small$\varepsilon_i\ee_i+\varepsilon_j\ee_j+\ee_k-\ee_l$]{}; 
\node (+-)	at 	(2,-2) [vx] [fill=black!50,label=below:\small$\varepsilon_i\ee_i+\varepsilon_j\ee_j-\ee_k+\ee_l$]{}; 
\node (--)	at 	(-2,-2) [vx] [fill=black!50,label=below:\small$\varepsilon_i\ee_i+\varepsilon_j\ee_j-\ee_k-\ee_l$]{}; 
\draw   (-1)--(2)--(1)--(-2)--(-1);
\draw  (+-)--(1)--(++)--(2)--(-+)--(-1)--(--)--(-2)--(+-);
\draw [thick] (++)--(+-)--(--)--(-+)--(++);
\end{tikzpicture}
\caption{Antiprism in the link of [$\varepsilon_i\ee_i,\varepsilon_j\ee_j]$ in $X^{(2)}$}
\label{antiprism}
\end{center}
\end{figure}
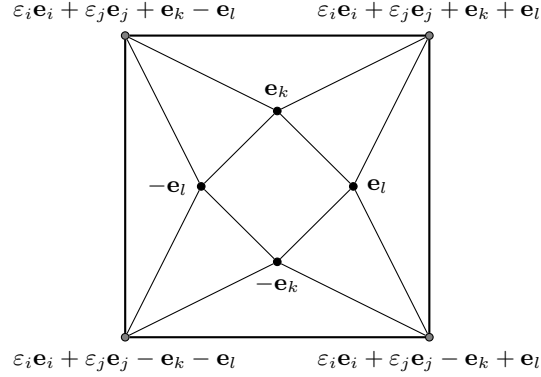

In either case, the boundary of the link of the edge $[\varepsilon_i\ee_i,\varepsilon_j\ee_j]$ is the boundary 
of the square
$$(\varepsilon_i\ee_i+\varepsilon_j\ee_j+\ee_k+\ee_l,\varepsilon_i\ee_i+\varepsilon_j\ee_j+\ee_k-\ee_l,
\varepsilon_i\ee_i+\varepsilon_j\ee_j-\ee_k-\ee_l,\varepsilon_i\ee_i+\varepsilon_j\ee_j-\ee_k+\ee_l).$$
We triangulate each of these squares by joining one pair of non-adjacent vertices by a diagonal.
Prima facie, we seem to have some amount of choice in this situation. All we have to ensure here 
is that if we introduce an edge $[q_{\varepsilon},q_{\varepsilon'}]$, then we also include the edge
$[-q_{\varepsilon},-q_{\varepsilon'}]$.

Recall that the $4$-dimensional hypercube $Q^4$ is bipartite, and any vertex-colouring partitions the  
vertices $\{q_{\varepsilon}|\varepsilon\in\{\pm1\}^4\}$ into two sets,
$$Q_e=\{q_{\varepsilon}|\varepsilon\in\{\pm1\}^4,\prod_{i=1}^4\varepsilon_i=+1\}\text{ and }
Q_o=\{q_{\varepsilon}|\varepsilon\in\{\pm1\}^4,\prod_{i=1}^4\varepsilon_i=-1\}.$$
Also note that if $q_{\varepsilon}$ and $-q_{\varepsilon}=q_{-\varepsilon}$ are always in the same 
block of the partition. Also any square in $Q^4$ contains exactly two vertices from $Q_o$ and two from $Q_e$. 

So we can triangulate each square with boundary 
$$(\varepsilon_i\ee_i+\varepsilon_j\ee_j+\ee_k+\ee_l,
\varepsilon_i\ee_i+\varepsilon_j\ee_j+\ee_k-\ee_l,
\varepsilon_i\ee_i+\varepsilon_j\ee_j-\ee_k-\ee_l,
\varepsilon_i\ee_i+\varepsilon_j\ee_j-\ee_k+\ee_l)$$
by joining its vertices in either of $Q_o$ or $Q_e$ by a diagonal. 

But, since the link of $[\ee_1,-\ee_1]$ is already a sphere, we need to triangulate the boundaries of the links
of $\pm\ee_1$ ``internally''. Additionally, this triangulation can not introduce a diagonal through the
interiors of either sphere, as any point $u\in Q$ at (Hamming) distance $3$ to a point $v$ is at distance $1$ to its
antipode $-v$. There is one way of triangulating an $8$-vertex $2$-sphere with the given partial $1$-skeleton without 
interior diagonals, i.e., the triangulation of the solid cube into five tetrahedra. So we triangulate each square in 
the link of $\pm\ee_1$ by joining the elements of say, $Q_o$, by edges. 

Now each element $u$ of $Q_o$ is joined to three other elements $v_1,v_2,v_3$ of $Q_o$ at distance $2$ from it. 
The other three elements of $Q_o$ at distance $2$ from $u$ are $-v_1,-v_2$, and $-v_3$. So none of the elements 
of $Q_o$ can be joined when triangulating the remaining squares. This forces us to triangulate the remaining 
squares by joining the vertices in $Q_e$ by an edge. This is possible, since for each vertex of $Q_e$, the three 
vertices at distance $2$ from it, which are across a square in the link of some $[\pm\ee_1,\pm\ee_i]$, have been 
ruled out in the previous step.

So for $\{i,j,k\}=\{2,3,4\}$, we replace the join of a line segment and the
boundary of a square in $X^{(2)}$, i.e,
\begin{multline*}
[\varepsilon_i\ee_i,\varepsilon_j\ee_j]*(\partial[
\varepsilon_i\ee_i+\varepsilon_j\ee_j+\ee_1+\varepsilon_i\varepsilon_j\ee_k,
\varepsilon_i\ee_i+\varepsilon_j\ee_j-\ee_1-\varepsilon_i\varepsilon_j\ee_k]*\\
\partial[
\varepsilon_i\ee_i+\varepsilon_j\ee_j-\ee_1+\varepsilon_i\varepsilon_j\ee_k,
\varepsilon_i\ee_i+\varepsilon_j\ee_j+\ee_1-\varepsilon_i\varepsilon_j\ee_k] ),
\end{multline*}
with the join of a new line segment with the boundary of another square, i.e,
\begin{multline*}
[\varepsilon_i\ee_i+\varepsilon_j\ee_j+\ee_1+\varepsilon_i\varepsilon_j\ee_k,
\varepsilon_i\ee_i+\varepsilon_j\ee_j-\ee_1-\varepsilon_i\varepsilon_j\ee_k]
*(\partial[\varepsilon_i\ee_i,\varepsilon_j\ee_j]*\\
\partial[
\varepsilon_i\ee_i+\varepsilon_j\ee_j-\ee_1+\varepsilon_i\varepsilon_j\ee_k,
\varepsilon_i\ee_i+\varepsilon_j\ee_j+\ee_1-\varepsilon_i\varepsilon_j\ee_k] )
\end{multline*}
So we are adding two $4$-simplices for every edge in $C^4$, $24\times2=48$ in
total. The $f$-vector of the boundary remains $[24,120,192,96]$, the same as
that of $X^{(2)}$. Call the boundary of the current complex $X^{(3)}$.

Now consider the link of a vertex from $C^4$ in $X^{(3)}$. The vertex-set of the
link of $\varepsilon_i\ee_i$ in $X^{(3)}$ is the set of vertices of the
hypercube $Q^4$ whose $i^{\mathrm{th}}$ co-ordinate is $\varepsilon_i$. These
vertices span a $3$-dimensional cube, and the triangles of the link of
$\varepsilon_i\ee_i$ in $X^{(3)}$ are ``halves'' of squares of $Q^4$.

Now consider opposite pairs of vertices of $C^4$. The boundaries of the links of
$+\ee_i$ and $-\ee_i$ are opposite (cubical) faces of $Q^4$. Moreover, the map
$\mathbf{x}\mapsto-\mathbf{x}$ swaps the triangulations of these cubes. So the
boundary of $X^{3}$ is an antipodal $S^3$. 

In order to triangulate the link of $\varepsilon_i\ee_i, 2\leq i\leq 4$,  we could
triangulate the interior of the $8$-vertex $2$-sphere (or triangulated cube)
described above by taking its cone at the point $-\varepsilon_i\ee_i$. In other
words, we take the join of the edge $[\pm\ee_i]$ with the $8$-vertex $S^2$
which is now the link of both $\ee_i$ and $-\ee_i$. 

This leaves the pair $\pm\ee_1$. The boundaries of the links of either vertex
is an $8$-vertex $S^2$, or the boundary of a cube triangulated by joining
``every other vertex by an edge''.
As mentioned above, we triangulate the links of each
vertex by splitting it into five tetrahedra, the vertices of four of
which have one element each of $Q_e$ and its three neighbouring elements of
$Q_e$. The vertices of the fifth are four elements of $Q_o$. 

Now we apply the map $x\mapsto -x$ on $Q$. This gives a triangulated $\RR P^4$. \qed
\end{cons}

We give one more way of triangulating $\RR P^4$ with $16$ vertices. Here we work with polyhedral complexes 
instead of the usual simplicial complexes. Again, the idea is to construct a $4$-dimensional ball with 
antipodal boundary, then to quotient via a restriction of the antipodal map. The description of the polyhedral 
complex used in this construction sacrifices rigour in the service of intuition. See~\cite[Appendix A]{SoB} 
for a more rigorous treatment.

\begin{cons}\label{c3}
Start with a (solid, $3$-dimensional) cube $Q^3$, embedded in $\RR ^4$ with vertices $(\pm1,\pm1,\pm1,0)$. Consider 
its suspension $SQ^3$ at the points $(0,0,0,\pm1)$. The boundary of this object is a $3$-dimensional polyhedral 
complex with ten vertices and $2\times 6=12$ (square-)pyramidal faces. The base of each of these pyramids is a 
face $S(\pm i)$ of the cube $Q^3$ given by $x_i=\pm1,x_4=0$, where $1\leq i\leq 3$. We avoid triangulating
the interior of $SQ^3$ for the time being.
 
First, we construct a ``dual'' cell complex $D$ outside $SQ^3$ by adding faces of increasing dimension, 
starting with points.

Corresponding to each face with base square $S(\pm i)$ and apex $(0,0,0,\varepsilon)$ of $SQ^3$, (where 
$\varepsilon\in\{\pm1\}$), take the point $3(\sigma_1,\sigma_2,\sigma_3,\varepsilon)$, where $\mathbf{\sigma}=
(\sigma_1,\sigma_2,\sigma_3)$ is the vector in $\RR^3$ taking value $\pm i$ at the $i$-th coordinate and $0$ 
elsewhere. So for example, the pyramid with apex $(0,0,0,1)$ and base $S(-2)$ gives the point $(0,-3,0,3)$. 
Corresponding to each of the $12$ faces of the boundary of the suspended cube, we get $12$ points.

Now join each pair of points in $D$ by an edge if the corresponding facets of $SQ^3$ intersect in a triangle or
square. This gives six edges corresponding to each square of $Q^3$. Additionally, $SQ^3$ has $2\times 12=24$
triangles corresponding to each point in $\{(0,0,0,\pm1)\}$ and each edge of $Q^3$. This gives $24$ more edges.

Next, consider the edges of $SQ^3$, of which there are $12+2\times 8$. The twelve edges of $Q^3$ give
twelve rectangles in $D$, of which the long pair of edges corresponds to the squares of $Q^3$ which intersect 
in this edge, while the short pair corresponds to the two triangles  of $SQ^3$ intersecting this edge. Each of
the remaining sixteen edges has as endpoints a vertex $v$ of $Q^3$ and a point in $\{(0,0,0,\pm1)\}$. Each such
edge gives a triangle of $D$, whose edges correspond to the three squares in $Q^3$ intersecting in $v$.

Now for each of the $8+2=10$ vertices of $SQ^3$, we add a polyhedron to $D$. The eight vertices of $Q^3$ 
give triangular prisms, whose faces are the three rectangles in $D$ corresponding to the three edges intersecting 
in this vertex, and the two triangles corresponding to the edges joining the vertex to each of 
$(0,0,0,\pm 1)$. Corresponding to either of the vertices $(0,0,0,\pm 1)$, we have an octahedron whose faces 
correspond to the $8$ edges of $SQ^3$ intersecting at the chosen vertex. 

We summarize the dualities described above in a table.

\vspace{12pt}
\begin{tabular}{cllc}
 \multicolumn{2}{c}{$SQ^3$}&\multicolumn{2}{c}{$D$}\\
Dim&Faces&Faces&Dim\\
$0$&$8+2$ points&$8$ prisms$+2$ octahedra&$3$\\
$1$&$12+2\times8$ edges&$12$ rectangles$+2\times8$ triangles&$2$\\
$2$&$6$ squares$+2\times12$ triangles&$6+2\times12$ edges&$1$\\
$3$&$2\times 6$ pyramids&$2\times6$ points&$0$\\
\end{tabular}
\vspace{12pt}

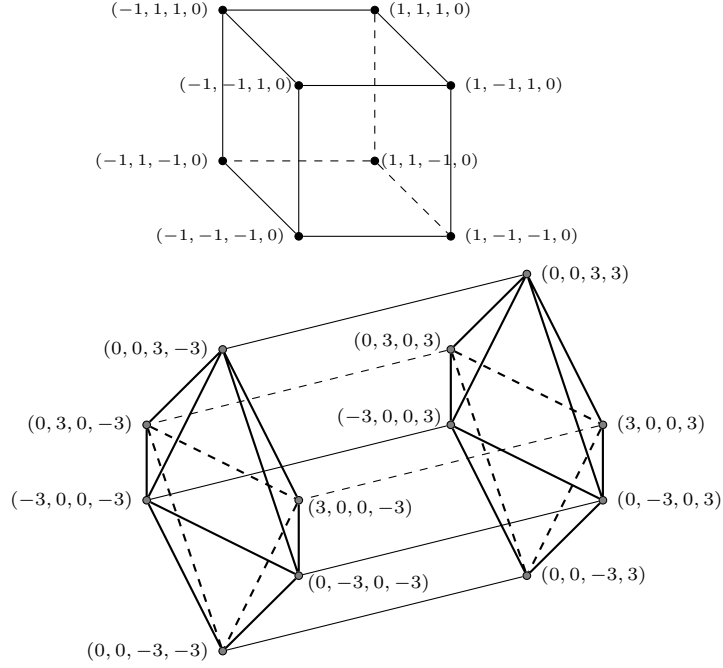
\begin{figure}%[b]
 \begin{center}
\begin{tikzpicture}[z=2cm,y=1cm,xslant=-1,
vx/.style={xshift=4cm,yshift=-4cm,circle,inner sep=0pt, minimum size=1mm,draw,fill=black!50},		      
vxx/.style={yshift=-5cm,circle,inner sep=0pt,minimum size=1mm,draw,fill=black!50},
v/.style={xshift=1.5cm, z=1cm, y=.5cm, circle,inner sep=0pt,minimum size=1mm,draw}]
\node (1) 	at 	(1,0,0) [vx] [label=right:{\scriptsize $(3,0,0,3)$}]{}; 
\node (-1) 	at 	(-1,0,0) [vx] [label={[xshift=2pt,yshift=3pt]left:{\scriptsize $(-3,0,0,3)$}}]{} ; 
\node (2)	at 	(0,1,0) [vx] [label={[xshift=2pt,yshift=3pt]left:{\scriptsize $(0,3,0,3)$}}]{} ; 
\node (-2)	at 	(0,-1,0) [vx] [label=right:{\scriptsize $(0,-3,0,3)$}]{} ; 
\node (3)	at	(0,0,1) [vx] [label=right:{\scriptsize $(0,0,3,3)$}]{} ; 
\node (-3)	at	(0,0,-1) [vx] [label=right:{\scriptsize $(0,0,-3,3)$}]{} ; 
\node (1') 	at 	(1,0,0) [vxx] [label= {[xshift=-2pt,yshift=-3pt]right:{\scriptsize $(3,0,0,-3)$}}]{}; 
\node (-1') 	at 	(-1,0,0) [vxx] [label=left:{\scriptsize $(-3,0,0,-3)$}]{} ; 
\node (2')	at 	(0,1,0) [vxx] [label=left:{\scriptsize $(0,3,0,-3)$}]{} ; 
\node (-2')	at 	(0,-1,0) [vxx] [label={[xshift=-2pt,yshift=-3pt]right:{\scriptsize $(0,-3,0,-3)$}}]{} ; 
\node (3')	at	(0,0,1) [vxx] [label=left:{\scriptsize $(0,0,3,-3)$}]{} ; 
\node (-3')	at	(0,0,-1) [vxx] [label=left:{\scriptsize $(0,0,-3,-3)$}]{} ; 
\draw [thick]  (-1)--(2)(1)--(-2)--(-1)--(3)--(1) (2)--(3)--(-2)--(-3)--(-1);
\draw [thick,dashed] (1)--(-3)--(2)--(1) (1')--(-3')--(2')--(1');
\draw [thick]  (-1')--(2')(1')--(-2')--(-1')--(3')--(1') (2')--(3')--(-2')--(-3')--(-1');
\draw (-1)--(-1') (-2)--(-2') (3)--(3') (-3)--(-3');
\draw [dashed] (1)--(1') (2)--(2');

\node (+++) 	at 	(1,1,1) [v] [fill,label=right:{\tiny $(1,1,1,0)$}]{}; 
\node (++-) 	at 	(1,1,-1) [v] [fill,label={[xshift=-3pt]right:{\tiny $(1,1,-1,0)$}}]{}; 
\node (+-+) 	at 	(1,-1,1) [v] [fill,label=right:{\tiny $(1,-1,1,0)$}]{}; 
\node (+--) 	at 	(1,-1,-1) [v] [fill,label=right:{\tiny $(1,-1,-1,0)$}]{}; 
\node (-++) 	at 	(-1,1,1) [v] [fill,label=left:{\tiny $(-1,1,1,0)$}]{}; 
\node (-+-) 	at 	(-1,1,-1) [v] [fill,label=left:{\tiny $(-1,1,-1,0)$}]{}; 
\node (--+) 	at 	(-1,-1,1) [v] [fill,label={[xshift=3pt]left:{\tiny $(-1,-1,1,0)$}}]{}; 
\node (---) 	at 	(-1,-1,-1) [v] [fill,label=left:{\tiny $(-1,-1,-1,0)$}]{}; 
\draw (---)--(+--)--(+-+)--(+++)--(-++)--(--+)--(---)--(-+-)--(-++)  (--+)--(+-+);
\draw[dashed](++-)--(-+-) (+++)--(++-)--(+--);
\end{tikzpicture}
 \end{center}
 \caption{Cube $Q^3$ and octahedral prism $D$}
 \label{octaprism}
\end{figure}

We can visualize the complex $D$ as a prismed octahedron, as in Figure~\ref{octaprism}. 

Now we can write down some of the $4$-simplices in our triangulation. Join each triangle in $SQ^3$ to its corresponding 
edge in $D$. This gives $24$ facets. Join each triangle in $D$ to its corresponding edge in $SQ^3$. 
This gives $16$ more facets.
  
Each vertex of $D$ is adjacent to five other vertices in $D$. In the $40$ simplices listed above, each vertex of $D$ is joined to 
four vertices of $Q^3$ and one vertex of $\{(0,0,0,\pm 1)\}$. Also, each of the vertices $(0,0,0,\pm1)$ is
joined to the six vertices of the octahedron corresponding to it in $D$. Each vertex of $Q^3$ is joined to 
three other vertices of $Q^3$, both of $(0,0,0,\pm1)$, and the six vertices of the triangular prism corresponding 
to it in $D$. Since the restricted antipodal map we wish to apply to the complex we are constructing takes $v\in D$ 
to $-v$, we can not add any more edges to our complex that contain a vertex of $D$. 

Now we consider the links of the $24$ triangles in $SQ^3$ in our complex so far. The triangle containing an edge whose 
$i,j^{\mathrm{th}}$ coordinates satisfy $x_i=\varepsilon_i, x_j=\varepsilon_j$, and the suspension point 
$(0,0,0,\varepsilon)$, is already joined to the (short) edge of $D$ whose endpoints correspond to the squares 
$x_i=\varepsilon_i,x_4=\varepsilon$, and $x_j=\varepsilon_j,x_4=\varepsilon$. Now we join each tetrahedron obtained by 
joining the above triangle to each of the vertices of the edge (in $D$) above, to a third point on the corresponding 
square in $Q^3$. Of the two remaining points on each square in $Q^3$, we choose the point the product of whose first 
three coordinates is $1$. For example, the triangle $[(1,1,-1,0),(1,-1,-1,0),(0,0,0,-1)]$ has as link $[(3,0,0,-3),
(0,0,-3,-3)]$. The corresponding simplices we add are the cones over tetrahedra
\begin{align*}
 [(1,1,-1,0)(1,-1,-1,0),(0,0,0,-1),(3,0,0,-3)]&\ast[(1,1,1,0)] \text{ and} \\
[(1,1,-1,0),(1,-1,-1,0),(0,0,0,-1),(0,0,-3,-3)]&\ast[(-1,1,-1,0)].
\end{align*}
Note that any such $4$-simplex can be obtained by starting with either of the two edges of $Q^3$ contained in it. 
So we get $2\times 12\times2/2=24$ simplices. 

The link of each of the $24$ triangles in $SQ^3$ is now a path of length $3$ with endpoints from the vertices of $Q^3$. 
Moreover, if the triangle $T$ is obtained by joining an edge $E$ of $Q^3$ with a point in $(0,0,0,\pm1)$, with $v\in E$ 
such that the product of the first three co-ordinates is $-1$, then the endpoints of the link of $T$ are the neighbours of
$v$ in the square in $Q^3$ containing $v$ but not $E$. For the next set of $4$-simplices, join the endpoints of the link
of each triangle is $SQ^3$ by an edge. For example, the triangle $[(1,1,-1,0),(1,-1,-1,0),(0,0,0,-1)]$ is joined to the edge
$[(1,1,1,0),(-1,1,-1,0)]$. So for each choice of $\{(0,0,0,\pm1)\}$ and each of the four vertices 
$v=(\varepsilon_1,\varepsilon_2,\varepsilon_3,0)$ of $Q^3$ such that $\varepsilon_1\varepsilon_2\varepsilon_3=-1$, 
we have a simplex whose vertices are $v$, its three neighbours in $Q^3$, and one of $(0,0,0,\pm 1)$. This gives eight more 
simplices. The link of every triangle in $SQ^3$ is now a circle. 

In constructing the previous two sets of simplices, we added one diagonal to each square face of $Q^3$. Recall that each 
square in $Q^3$ given by the equation $x_i=\varepsilon_i$ corresponds to the edge of $D$ spanned by 
$3(\sigma_1,\sigma_2,\sigma_3,\pm1)$, where $\sigma_j=\varepsilon_i$ if $i=j$ and $0$ otherwise. Now consider a
triangle in our complex with vertices in this square, say $[v_0,v_1,v_2]$, and let the product of the first three 
coordinates of $v_0$ be $-1$. The link of this triangle has four edges. Let $v_0'$ denote the third vertex adjacent to 
$v_0$ in $Q^3$. In the link of this triangle, $v_0'$ is joined to $(0,0,0,\pm1)$, and the latter vertices are respectively
joined to $3(\sigma_1,\sigma_2,\sigma_3,\pm1)$. We join the triangle $[v_0,v_1,v_2]$ to the line segment 
$[3(\sigma_1,\sigma_2,\sigma_3,1),3(\sigma_1,\sigma_2,\sigma_3,-1)]$. The twelve triangles with vertices on a
square in $Q^3$ gives one simplex each, so we get twelve new simplices. 

Also, in our last but one set of simplices, we introduced four new triangles, each consisting of three vertices of $Q^3$, 
such that the first three coordinates of each have product $1$. The link of each such triangle consists of two edges,
where the common neighbour of the three vertices is joined to each of $(0,0,0,\pm1)$. So we have four triangles forming 
the boundary of a tetrahedron, the boundaries of the links of each being the vertices $(0,0,0,\pm1)$. We add two new 
simplices, by taking the tetrahedron consisting of these four vertices of $Q^3$ and joining it to each of the points 
$(0,0,0,\pm1)$.

So now the links of all triangles with vertices from $SQ^3$ are circles. We consider the links of edges in $SQ^3$. 

The link of an edge of the form $[(\varepsilon_1,\varepsilon_2,\varepsilon_3,0),(0,0,0,\varepsilon)]$, where 
$\varepsilon_1\varepsilon_2\varepsilon_3=-1$, is an octahedron. One of the faces of this octahedron consists of the 
three neighbours of $(\varepsilon_1,\varepsilon_2,\varepsilon_3,0)$ in $Q^3$ and its opposite face in this octahedron
is the triangle corresponding to the chosen edge in $D$. 

The link of an edge of the form $[(\varepsilon_1,\varepsilon_2,\varepsilon_3,0),(0,0,0,\varepsilon)]$, where 
$\varepsilon_1\varepsilon_2\varepsilon_3=1$, is a $9$-vertex $S^2$ with six vertices from $Q^3$ and the 
remaining three vertices from its corresponding triangle in $D$. 

The boundary of the link of an edge in $Q^3$ is the boundary of its corresponding square in $D$.

Note that the boundary of this complex is now a triangulated $S^3$ with $f$-vector $[22,102,160,80]$.

Now consider each edge $[v_1,v_2]$ in $Q^3$ and its opposite edge in $Q^3$, $[-v_1,-v_2]$. If the boundary of
the links of the first edge is $[w_1,w_2],[w_2,w_3],[w_3,w_4],[w_4,w_1]$, then the link of the opposite edge
is $[-w_1,-w_2],[-w_2,-w_3],[-w_3,-w_4]$, $[-w_4,-w_1].$ 

Now we apply the antipodal map $v\mapsto-v$ on the vertices of $D$, the two above squares will be identified. 
So will the two octahedra $O^+$ and $O^-$ which are the boundaries of the links of $\pm(0,0,0,1)$ respectively.

We close the links of the edges $\pm[v_1,v_2]$ by joining each of the tetrahedra containing $[v_1,v_2]$ to that 
vertex of $[-v_1,-v_2]$, the product of whose first three coordinates is $-1$. 

We close the boundary by adding the simplices obtained by joining each of the triangles $[(1,1,1,0)$, $(1,1,-1,0)$,
$(-1,-1,-1,0)]$ and $[(-1,-1,-1,0)$, $(-1,-1,1,0)$, $(1,1,-1,0)]$ with each of the edges in $S$.
This gives $2\times6\times4=48$ simplices.

Now we have joined each vertex $v$ in $Q^3$ to its opposite vertex $-v$. The link of the edge $[v,-v]$ consists of 
twelve triangles. Suppose $v=(\varepsilon_1,\varepsilon_2,\varepsilon_3,0)$ where 
$\varepsilon_1\varepsilon_2\varepsilon_3=1$, and let $\overline{P_v}$ be the image under the antipodal map on $D$
of the prism(s) corresponding to $v$ (and $-v$) in $D$. Then the faces of the link of $[v,-v]$ are the three square faces
of $\overline{P_v}$, each subdivided by the corresponding neighbour of $v$. The boundary of this complex is the set of 
edges of two disjoint triangles in the image of $O^+$ (and $O^-$) under the map $x\mapsto-x$. We add the joins of
$[v,-v]$ with each of these triangles. This gives $4\times 2=8$ simplices.

We close the boundaries of $\pm(0,0,0,1)$ by joining each of the faces of the octahedron to the edge 
$[(0,0,0,1),(0,0,0,-1)]$. This gives eight more simplices. 

This gives a $16$-vertex triangulation of $\RR P^4$ with $150$ simplices. \qed
\end{cons}

\section{Similarities, differences, and concluding remarks}

The starting objects of our two latter constructions, such as the hyperoctahedron, the $4$-cube, the suspended 
cube, and octahedral prism suggest that automorphism groups of the triangulations we constructed are very close to
$C_2\times S_4$. In fact the simplices we add in each step of each construction are typically orbits under this group.
But on explicit computation, the automorphism groups of both complexes turn out to be isomorphic to $S_6$, acting 
on $10+6$ vertices. 

In Construction~\ref{c2}, the vertex-orbit of $S_6$ of size $6$ consists of the two points of the hyperoctahedron 
$C^4$ used to triangulate it internally, (namely $\pm\mathbf{e}_1$), and the vertices of $Q_e$.

In Construction~\ref{c3}, the smaller $S_6$ orbit consists of the two suspension points $(0,0,0,\pm1)$, 
and the four vertices of the cube the product of whose first three coordinates is $-1$. 

One construction starts with a suspended octahedron (hyperoctahedron) on the inside and a cubical prism (hypercube) on 
the outside of our $4$-dimensional ball. The other starts with a suspended cube on the inside and an octahedral prism 
on the outside. This gives us a way of visualizing either construction as the other one ``turned inside-out''. 

Also recall that the $C_2\times S_4$ is the stabilizer of a $2$-subset in $S_6$. If we consider any pair of elements
of the orbit $\mathcal{O}_6$ of size $6$ in either construction, we find that the link of the edge joining them is an
octahedron consisting of six points of the longer orbit $\mathcal{O}_{10}$, and that the intersection of their 
links is a solid cube triangulated with five tetrahedra, where the vertices of the inner tetrahedron are the remaining 
vertices of $\mathcal{O}_{10}$ and the other four vertices are the remaining vertices of $\mathcal{O}_6$. This also 
gives a correspondence with the vertex set of the triangulated $\RR P^4$ of Construction~\ref{c1}, which induces 
simplicial isomorphisms between all three complexes. This gives weight to our belief that the following is true.

\begin{conj}\cite{D}
 $\RR P^4_{16}$ is the unique $16$-vertex triangulation of $\RR P^4$.
\end{conj}

Their close connection not withstanding, Constructions~\ref{c2}~and~\ref{c3} offer different perspectives on how to construct 
triangulations of $\RR P^n$ for other values of $n$. The boundaries of the very final $4$-balls we construct before applying the 
antipodal map are antipodal $3$-spheres. In the first case, the $3$-sphere has $24$ vertices, and the quotient only
gives a $12$-vertex $\RR P^3$, whereas in the second construction, we get the same $3$-sphere constructed by Walkup 
as the double cover of his $\RR P^3_{11}$.

Also note that Walkup's $\RR P^3_{11}$ also fits into the paradigm outlined in Construction~\ref{c3}. 
Recall that the initial object of our construction was an octahedron, surrounded by a cube. Observe that the
octahedron and cube are respectively a suspension of a square and a prismed (solid) square. 

The two constructions in Section~\ref{sec:fc} point to different possible generalizations. Construction~\ref{c2} 
suggests constructions using a hyperoctahedron placed within a hypercube to obtain an $\RR P^n$ on $2^{(n-1)}+2n$
or $2^{(n-1)}+2n+1$ vertices. This author has tested this approach to construct an $\RR P^5$ on $2^{(5-1)}+2\times
5+1=27$ vertices and an $\RR P^6$ on $2^5+12+1=45$ vertices. 

Construction~\ref{c3} suggests the possibility of taking a mixture of suspended and prismed polyhedra in higher dimensions to 
lower the number of vertices needed even further. If this can be realized by starting with a double suspension of the $3$-cube, 
then it may be possible to triangulate $\RR P^5$ with only $8+4+\frac{24}{2}=24$ vertices. F.H. Lutz has discovered an 
$S_4$-invariant $24$-vertex triangulation of $\RR P^5$, the facet-list of which is available at~\cite{Lw}. This author has not 
investigated the possibility that this complex is isomorphic to the one we postulate. The {\tt{BISTELLAR}} program used to 
discover this triangulation is available through the {\tt{GAP}} package {\tt{simpcomp}}~\cite{sc}, along with a library of known 
triangulations and other tools for calculating with complexes.

Our observation that double covers of $\RR P^{n-1}$ appear as boundaries of $n$-balls in the constructions of $\RR P^n$ 
suggests the following question. Given an antipodal $n$-sphere, is it possible to ``thicken'' it to an $n+1$-ball while 
preserving a suitable restriction of the antipodal map, then glue the quotiented boundary to itself to get a triangulated 
$\RR P^{n+1}$? If the answer is yes, it suggest the possibility of inductively triangulating $\RR P^n$. In the case of 
the double cover of $\RR P^4_{16}$ in Construction~\ref{c1}, it is easy to construct a $B^5$ on $32$ vertices such that 
the antipodal map applies to the orbit of size $20$. But the question of whether there exists a gluing of the quotiented 
boundary which produces a triangulated $\RR P^5$ is harder to answer. If such a triangulation exists, it
would have $12+\frac{20}{2}=22$ vertices, which would make it vertex-minimal by the theorem of Arnoux and Marin.

\subsection*{Acknowledgments}
This work is based on part of my PhD thesis, written at Maynooth University, supported by a John and Pat Hume Research 
Scholarship. I would like to thank Pat McCarthy, my thesis advisor, for many fruitful discussions, which led to some of 
the results in this paper. Section~\ref{sec:fc} is essentially joint work with him. I would also like to thank Mikhail 
Klin for pointing out the relevant history of the construction of $W_{22}$, and Stefan Bechtluft-Sachs and Eran Nevo, 
for their helpful comments which improved the presentation of the paper.

%Conflict of Interest: The author declares that she has no conflict of interest.

\bibliographystyle{plain}
\bibliography{rp4.bib}

\end{document}